\documentstyle [12pt,epsf,leqno] {article}                                       
\def\RR{I\kern-0.35em R\kern0.2em \kern-0.2em}
\def\BBN{I\kern-0.35em  N\kern0.2em \kern-0.2em}
  
\newcommand{\rdr}{\, r\, dr}
\newcommand{\EE}{{\cal E}_{\kappa,g}}
\newcommand{\EEhk}{{\cal E}_{\infty,g}}
\newcommand{\II}{I_{\kappa,g}}
\newcommand{\IIhk}{I_{\infty,g}}

\newcommand{\weak}{\rightharpoonup}

\newcommand{\eps}{\varepsilon}
\renewcommand{\epsilon}{\varepsilon}
\renewcommand{\phi}{\varphi}

\newcommand{\QED}{\newline $\diamondsuit$}

\newcommand{\CC}{C_0^\infty}
\newcommand{\Pf}{{\bf Proof:}\  }

\newcommand{\GL}{(GL)${}_{\kappa,g}$}
\newcommand{\GLhk}{(GL)${}_{\infty,g}$}
\newcommand{\fk}{\tilde f_\kappa}
\newcommand{\finf}{\tilde f_\infty}
\newcommand{\Sk}{\tilde S_\kappa}
\newcommand{\fhk}{\hat f_\kappa}
\newcommand{\Shk}{\hat S_\kappa}
\newcommand{\mhk}{\hat m_\kappa}
\newcommand{\uk}{u_\kappa}
\newcommand{\k}{\kappa}
\newcommand{\FF}{{\cal F}}
\newcommand{\LL}{{\cal L}}
\newcommand{\codim}{{\rm codim}\,}
\newcommand{\Ran}{{\rm Ran}\,}
\newcommand{\ep}{\epsilon}
\newcommand{\rhe}{{\rho_\epsilon}}

\newcommand{\be}{\begin{equation}}
\newcommand{\ee}{\end{equation}}
\newcommand{\bea}{\begin{eqnarray}}
\newcommand{\eea}{\end{eqnarray}}
\newcommand{\beann}{\begin{eqnarray*}}
\newcommand{\eeann}{\end{eqnarray*}}
\newcommand{\nnn}{\nonumber}

\begin{document}

\newtheorem{th}{Theorem}[section]
\newtheorem{pr}[th]{Proposition}
\newtheorem{lem}[th]{Lemma}
\newtheorem{de}[th]{Definition}
\newtheorem{re}[th]{Remark}
\newtheorem{co}[th]{Corollary}

\baselineskip=17pt
\begin{titlepage}
\title{  Vortex Structures for an $SO(5)$ Model of
High-$T_C$ Superconductivity and Antiferromagnetism \\  }
\author{ {\Large Stan Alama} \thanks{Supported by an
NSERC (Canada) Research grant. e-mail: {\tt alama@mcmaster.ca}} 
\and
{\Large Lia Bronsard} \thanks{Supported by an
NSERC (Canada) Research grant. e-mail: {\tt bronsard@math.mcmaster.ca}}  
\and {\Large Tiziana Giorgi} \thanks{e-mail: {\tt giorgi@math.mcmaster.ca}}  \\
{\small McMaster Univ., Dept. of Math. \& Stat., Hamilton Ont.,
L8S 4K1 Canada}\normalsize }

\thispagestyle{empty}
\maketitle

\begin{abstract}
We study the structure of symmetric vortices in a
Ginzburg--Landau model based on S.\ C.\ Zhang's $SO(5)$
theory of high temperature superconductivity and antiferromagnetism.
We consider both a full Ginzburg--Landau theory (with
Ginzburg--Landau scaling parameter $\kappa<\infty$) and a $\kappa\to
\infty$ limiting model.
In all cases we find that the
usual superconducting vortices (with normal phase in the
central core region) become unstable (not energy minimizing)
when the chemical potential crosses a threshold level, giving
rise to a new vortex profile with antiferromagnetic ordering
in the core region.  We show that this phase transition in
the cores is due to a bifurcation from a simple eigenvalue
of the linearized equations.  In the limiting large $\kappa$
model we prove that the antiferromagnetic core solutions
are always nondegenerate local energy minimizers and
prove an exact multiplicity result for physically relevent solutions.
\end{abstract}

\end{titlepage}
\newpage

\section{Introduction}

In 1986 Bednorz and M\"uller
  announced their discovery of
high critical-temperature ($T_C$)
superconductors, and promptly received the 1987 Nobel Prize for
their efforts.  This discovery has led 
to a new flowering of superconductivity theory, since the high
temperature phenomenon cannot be explained by the
accepted models  for conventional superconductors. 
In particular, many physicists
have come to the conclusion that the microscopic BCS theory does not
correctly describe the interactions which produce superconductivity at
high temperatures. At the present time, there are several competing theories
which attempt to explain these interactions. One theory is based on the
observation that high-$T_C$ compounds also exhibit an ordered phase called
{\it antiferromagnetism} when physical parameters (such as
temperature, chemical potential or ``doping'', and magnetic field) are varied.
Antiferromagetism (abbreviated AF) is an insulating phase of matter in which
electron spins orient themselves in the direction opposite to their
nearest neighbors.   The coexistence of these two phases (AF and SC)
in the phase diagram of the high-$T_C$ compounds has led to the
speculation that high temperature superconductivity and 
antiferromagnetism could be
explained by the same type of interaction.  

Following in this
direction, Shou-Cheng Zhang \cite{Z}  proposed a quantum
statistical mechanics model 
which incorporates AF and high temperature superconductivity (SC). The model
is based on a broken $SO(5)$ symmetry tying the complex 
order parameter of superconductivity to the N\'eel vector which
describes antiferromagnetism.
The interactions between the SC
and AF order parameters in this model should have some effect
on the familiar constructions from conventional superconductivity theory.
In a recent paper Arovas, Berlinsky, Kallin, \& Zhang \cite{ABKZ} 
introduced a phenomenological Ginzburg--Landau model
based on the $SO(5)$ theory, and studied isolated vortex solutions
in the plane.  Recall that in a conventional superconductor the magnetic
field is expelled from the superconducting bulk, and only
penetrates in thin tubes (the vortices) where superconductivity
is supressed. Hence, in the conventional theory
 the magnetic field is constrained
to a small {\it core} of normal (non-SC) phase.
Using a simplified model Arovas {\it et al} predicted a
new kind of vortex structure in the $SO(5)$ model:  vortices
with antiferromagnetic cores, which should be observed
for small values of the chemical potential.  They also predicted 
that (as the chemical potential is gradually decreased) the transition from
normal core to AF core vortices occurs in a discontinuous fashion.  In
other words,  AF cores should be produced via a  
{\it first order} phase transition.

In this paper we rigorously analyse vortex cores in the full
$SO(5)$ Ginzburg--Landau model and in an ``extreme type II''
limiting model (also called ``high kappa model'') to understand the nature of
the transition between normal core and AF core solutions.
For both models we show that the vortex solutions with
normal cores become unstable (within the class of radial
functions-- see (\ref{radial}) below,) and vortices with AF cores
are produced by bifurcation from the normal core solutions.  In the extreme
type II model we prove that the transition is continuous
(ie, {\it second order}), contrary to the prediction of \cite{ABKZ}
(see Figure 1.)  Furthermore,
we show that for each value of the chemical potential there exists
a unique stable vortex profile  (see Theorem~\ref{unique}.)
\medskip
 
The full $SO(5)$ Ginzburg--Landau free energy is written in 
terms of the SC order parameter $\psi\in {\bf C}$ and
the AF order parameter (N\'eel vector) $\vec m=(m_1,m_2,m_3)$.
In non-dimensional form, the free energy is:
$$  
{\cal F} =
\frac12 \int_\Omega \left\{  {\kappa^2\over 2}(1-  |\psi|^2 -|{\vec m}|^2)^2 
    + g\kappa^2 |\vec m|^2
   + |(\frac{1}{i}\nabla-{\vec A})\psi|^2 + 
         |\nabla {\vec m}|^2 + |\nabla\times {\vec A}|^2
    \right\} \, dx.          
$$
(We refer to the paper by Alama, Berlinsky, Bronsard \&
Giorgi \cite{ABBG} where the free energy is written in
dimensional form.)
In these variables, the penetration depth $\lambda=1$, and the Ginzburg--Landau
parameter
$\kappa$ is the reciprocal of the correlation length $\xi$.
The parameter $g$ measures the strength of doping 
(chemical potential) of the material.
It is this term which breaks the $SO(5)$ symmetry of the potential term.
We take $g>0$: with this assumption superconductivity is preferred in the bulk
of the sample.

To study isolated vortex solutions in the plane 
$\Omega=\RR^2$ we seek critical points
of $\cal F$ of the form
\be\label{radial}
  \psi = f(r)e^{id\theta}, \qquad
     {\vec A}= S(r) \, \left( {-y\over r^2}, {x\over r^2} \right), 
      \qquad {\vec m}=m(r){\vec  m_0} 
\ee
where ${\vec  m_0}$ a fixed unit vector, and $d\in {\bf Z}\setminus\{0\}$ represents
the degree of the vortex.  As for conventional SC vortices, we expect
that only the solutions with $d=\pm 1$ will be energy minimizers
(see Gustafson \cite{Gustafson}, Ovchinnikov \& Sigal \cite{OS}.)
Critical points of $\cal{F}$ with this ansatz solve the system of equations
$$
\left\{  \begin{array}{c}
 -f'' - {1\over r}f' + {(d-S)^2\over r^2}f = \kappa^2 (1-f^2-m^2)f, \\  \\
  -S'' + {1\over r} S' = (d-S)f^2,  \\ \\
-m'' -{1\over r} m' + \kappa^2 g m = \kappa^2 (1-f^2-m^2)m,
\end{array} \right. \leqno(GL){}_{\kappa,g}
$$
with $f(r)\ge 0$, $f(r), S(r)\to 0$ as $r\to 0$, and $f(r)\to 1$; $S(r)\to d$ as
$r\to\infty$; and $m'(0)=0$, $m(r)\to 0$ as $r\to\infty$. 

In addition, we study the following ``extreme Type II'' model, 
$$
\left\{  \begin{array}{c}
 -f'' - {1\over r}f' + {d^2\over r^2}f = (1-f^2-m^2)f, \\ \\
-m'' -{1\over r} m' +  g m = (1-f^2-m^2)m.
\end{array} \right. \leqno(GL){}_{\infty,g}
$$
The system \GLhk\ is obtained in the limit $\kappa\to\infty$
after rescaling solutions to \GL\ by the correlation length
$\xi=1/\kappa$.  For high $T_C$ superconductors $\kappa$ is
very large, and hence the vortex cores are very narrow compared
to the penetration depth, which measures the length scale
for magnetic fields.  By rescaling we capture the structure
of the vortex cores and decouple the magnetic field, which
lives on a much larger length scale.  Indeed, the calculations
which led Arovas {\it et al} \cite{ABKZ} to predict AF
vortex cores are mostly based on \GLhk\ and its associated
free energy functional.

\medskip

We observe that when the AF order parameter $m=0$
the two systems \GL\ and \GLhk\ reduce to the familiar
Ginzburg--Landau vortex equations, well studied in the
mathematical literature (see Plohr \cite{Plohr}, Berger and Chen \cite{BC},
Chen, Elliot, \& Qi \cite{CEQ}, Brezis, Merle, \& Rivi\`ere \cite{BMR},
Ovchinnikov \& Sigal
\cite{OS},  for example.)  We call these the {\it normal core} solutions.
In a previous paper \cite{ABG} we have proven that
when $\kappa^2\ge 2d^2$ there is a unique normal core
solution, which is a non-degenerate minimizer of the 
appropriate free energy functional.  This characterization
will be essential for our analysis of the normal-to-AF core
transition.

\medskip

We now discuss our results.
We define a reduced energy functional defined for functions
satisfying the symmetric vortex ansatz (\ref{radial}), as
well as appropriate function spaces in which that functional
is smooth.
We find that for every
$\kappa$ (including the extreme type II model) there
exists $g^*_\kappa>0$ such that the conventional
normal core vortex solutions of \GL\ (and \GLhk)
are strict local  minimizers of the reduced energy  for $g>g^*_\kappa$, but are
not local minimizers when $0<g<g^*_\kappa$.  In particular,
energy minimizers must have AF order in the vortex core
for $0<g<g^*_\kappa$.  When $\kappa^2\ge 2d^2$ we show that the AF core
solutions bifurcate from the normal core solution at
a simple eigenvalue of the linearized system \GL\
(or \GLhk.)  The bifurcating solutions remain bounded for
$g>0$ and lose compactness as $g\to 0+$ with $f\to 0$ and $m\to 1$.

 For the limiting problem
\GLhk\ we obtain a complete picture of the phase
transition to AF cores.  This is because all
AF core vortex solutions
are non-degenerate minima of the reduced energy. 
(See Theorem~\ref{nondeg}.)
Stable (locally minimizing)
solutions with $m(r)>0$ bifurcate from $m=0$ at
$g=g^*_\infty$ to values $g<g^*_\infty$.
 Moreover, for each $g< g^*_\infty$ there exists exactly
one solution with $m(r)>0$.  

In the language of physics, our results indicate a {\it second
order} (or continuous) phase transition between
normal and AF vortex cores in \GLhk.  This information concerning the
nature of the transition was not derived in the paper
by Arovas {\it et al} \cite{ABKZ}, and hence the result is
new to the physics literature as well.  
For \GL\ 
Alama, Berlinsky, Bronsard, \& Giorgi \cite{ABBG} present
numerical simulations (based on gradient flow for a finite elements
approximation of the free energy) which suggest
that the transition is also second order for $\kappa<\infty$.  (See Figure 1.)
However we were not able to extend the arguments used in studying
the bifurcation curves of \GLhk\ to the more complicated system
\GL.  See Remark~\ref{bifdir} for further discussion.

\medskip

Here is an outline of the content of the paper.
In the second section we introduce the reduced energy and function
spaces, we treat briefly  the questions of existence, regularity, and decay of solutions,
and we present properties of physically relevant (``admissible'') solutions.
We also prove the monotonicity of the solution profiles $(f,S,m)$ under
the hypothesis that the solution is a local reduced energy minimizer.
This result (Theorem~\ref{monotonicity}) is done in the spirit
of the weak maximum principle (see Theorem~8.1 of \cite{GT}.)

Section 3 contains the proof that all solutions of \GLhk\ with
$m>0$ represent non-degenerate local minima of the reduced
energy.  This result is the key to understanding the bifucation
diagram for \GLhk.  The bifurcation analysis itself occupies Section 4.

The last two sections contain the {\it a priori} estimates used in
rigorously passing to the
limit $\kappa\to\infty$ and in
studying the global behavior of bifurcating continua.
In both cases, we require estimates on solutions which are
energy-independent.  For the limit $\kappa\to\infty$ this
is because the reduced energy of minimizers behaves like
$\log\kappa$, and in studying global bifurcation we require
estimates valid for any physically relevant solution (whether
it is energy minimizing or not.)  
The starting point for these estimates is a Pohozaev type identity 
(see Proposition~\ref{pohozaev}.)
The proof of convergence to \GLhk\ as $\kappa\to\infty$ is
presented in Section 5;  other {\it a priori} estimates are derived
in Section 6.

\medskip

We wish to thank our colleague John Berlinsky for introducing
us to the $SO(5)$ model, and for 
his great patience in explaining physics to we
mathematicians.  We are also obliged to the Brockhouse Institute
for Materials Research for supporting a workshop which brought
together physicists and mathematicians to discuss 
issues in superconductivity.

\begin{figure} 
\begin{center}
\leavevmode
\hbox{%
\epsfxsize=5.5in
\epsffile{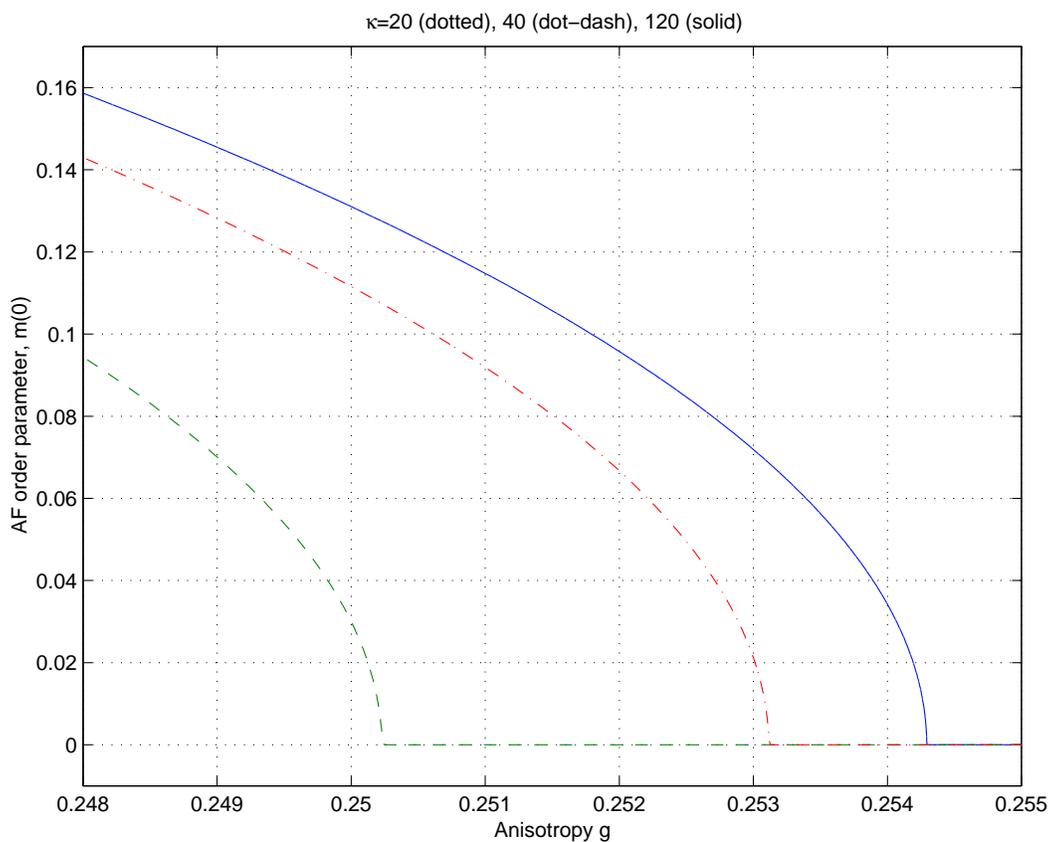}}
\small
\caption{A numerical bifurcation curve, $m(0)$ vs. $g$ for values $\kappa=20, 40, 120$
and $d=1$, indicates a second-order transition to AF cores in model
\GL. For the high-kappa model \GLhk\ we prove that the above image
correctly depicts the solution set (see Theorem~\ref{unique}.)
Numerical simulations indicate that the bifurcation
occurs at $g_\infty^*\simeq 0.2545$
\cite{ABBG}.}
\label{fig1}
\end{center}
\end{figure}

\section{Solutions of the Ginzburg--Landau system}

\setcounter{equation}{0}
\setcounter{th}{0}

\subsection{Preliminaries}

Here and in the rest of the paper, we fix the value of $d\in {\bf Z}\setminus\{0\}$.
In this section $\kappa\in\RR$ is fixed.  Note that without loss of generality
we may take $d>0$, since the free energy and the corresponding
Euler--Lagrange equations are invariant under
the transformation $(\psi, {\bf A}, {\bf m})\to(\bar\psi, -{\bf A}, {\bf m})$.

Following our previous work \cite{ABG} on symmetric
vortices, we define a function space  for which the free energy will be a smooth
functional.  First we fix some notation:  we denote
by $L^p_r$, $H$ the Lebesgue and Sobolev spaces (respectively)
 of radially symmetric functions in $\RR^2$, that is,
\beann
L^p_r &=& \{  u(r):
   \int_0^\infty |u(r)|^p\rdr <\infty \}, \quad (p<\infty), \\
H:=H^1_r &=& \{u(r): \int_0^\infty \left[ (u'(r))^2 + (u(r))^2\right]\rdr <\infty\},
\eeann
and analogously for $L^\infty_r$.  We also denote
$\int u(r)\rdr = \int_0^\infty u(r)\rdr$.

Define the Hilbert space
$$    X =\{ u\in H: \
   \int {u^2\over r^2}\rdr <\infty \},  $$
with norm
$$   \|u\|_X = \sqrt{  \int \left[
       (u'(r))^2 +u^2 +{u^2\over r^2}\right] \rdr. }
$$
The following density and imbedding properties for 
the space $X$ are proven in \cite{ABG}:
\begin{lem}\label{embedding}
\begin{enumerate}
\item  $X$ is compactly embedded in $L^p_r$ for each
$p \in (2,\infty)$.
\item  $X$ is compactly embedded in $L^2_{r,loc}$.
\item   For every $u\in X$,
$$  \|u\|^2_\infty  \le \int \left[
     (u')^2 + {u^2\over r^2} \right] \rdr.
$$
In particular, $X$ embeds continuously into $L^\infty_r$.
\item $\CC((0,\infty))$ is dense in $X$.
\end{enumerate}\end{lem}
We note that the compactness of the embedding of $H$
into $L^p_{r,loc}$ $(1\le p <\infty)$ is just the classical
Rellich-Kondrachov Theorem, and the compact embedding
of $H$ into $L^p_r$ for $2<p<\infty$ is due to Strauss
\cite{Str}.

\subsection{Energy}

We now define our energy functionals, using the space $X$ defined
above.  To keep the appropriate boundary condition at
infinity we fix any function $\eta\in C^\infty([0,\infty))$ with
$\eta(r)=0$ for $0\le r\le 1$, $\eta(r)=1$ for all $r\ge 2$, and $0<\eta<1$.
Then set $f_0=\eta$, $S_0=d\,\eta$, and seek solutions $(f,S,m)$
of \GL\ with $f=f_0+u$, $S=S_0+ rv$, $u,v\in X$, $m\in H$.
(Later we will see that this choice poses no restriction on solutions
which are physically relevant.)  We denote by $Y_0=X\times X\times H$, and by
$Y$ the affine space
$$  Y=\{ (f,S,m): \ f=f_0+u, \ S=S_0+rv,\  u,v\in X, \  m\in H\} = Y_0 + (f_0,S_0,0).  $$

For $(f,S,m)\in Y$ we define
\bea\label{energy}
&&  \EE (f,S,m) = \\
\nnn
&&  \quad
\frac 12 \int\left\{  (f')^2 + \left[ {S'\over r}\right]^2 + (m')^2 + \kappa^2 g m^2
     + {(d-S)^2 f^2\over r^2} + {\kappa^2\over 2}(1-f^2-m^2)^2 \right\}\rdr
\eea
and the functional $I_{k,g}: Y_0\to \RR$ by
$$     I_{k,g}(u,v,m) = \EE (f_0+u,S_0+rv,m)-\EE(f_0,S_0,0).
$$
Throughout the paper we will take advantage of these two representations of our
spaces and energies, and use the formulation which is more
convenient at the given moment.

Defining an energy functional for the limiting problem \GLhk\
is trickier, since the naive choice for the energy
(namely (\ref{energy}) with  $S=0$ and $\kappa=1$)
would be infinite for all $f$ satisfying the desired boundary
condition at $r=\infty$.  Our solution is to subtract off
the offending term from the energy density.
Let $\finf$ be the (unique) positive solution to the
high kappa vortex equation,
$$   -\finf'' -{1\over r}\finf' + {d^2\over r^2}\finf = 
(1-\finf^2)\finf
$$
with $\finf (0)=0$, $\finf (r)\to 1$ as $r\to\infty$.
The uniqueness of $\finf$ was established 
by Chen, Elliot \& Qi \cite{CEQ}.
The estimates in \cite{CEQ} ensure that $\finf$ is
smooth, $\finf(r)\sim r^d$ near $r=0$, and
$(1-\finf)\in H$.

We define the appropriate spaces for the free energy 
$\EEhk$ based on $\finf$:  let $Z_0=X\times H$ and
$$  Z= \{(f,m): \  f=\finf + u, \ u\in X, \ m\in H\}= Z_0 + (\finf,0).  $$
Then the energy for the high kappa model is:
\bea\label{hkenergy}
&&  \EEhk (f,m) = \\
\nnn
&& \quad
\frac 12\int\left\{(f')^2 +  (m')^2 +  g m^2
     + {d^2 \over r^2}[f^2-\finf^2] + {1\over 2}(1-f^2-m^2)^2 \right\}\rdr 
\eea
If we write $f=\finf+u$, we reduce to the equivalent functional
\bea
\label{hkI}
&& I_{\infty,g}(u,m)=   \EEhk(\finf + u, m)-\EEhk(\finf,0) \\
\nnn
&&  \quad
=\frac 12 \int \left\{  (u')^2 + {d^2\over r^2}u^2 + (m')^2 + gm^2  \right. \\
\nnn
&& \qquad  \left.
        + {1\over 2}(1-(\finf+u)^2-m^2)^2 
      - {1\over 2} (1-\finf^2)^2 + 2(1-\finf^2)\finf u \right\}\rdr.
\eea
By a direct expansion of the energy in powers of $u,v,m$ we see that
$I_{k,g}: Y_0\to \RR$ and $I_{\infty,g}: Z_0\to\RR$ are
smooth ($C^\infty$) functionals. 

When $g>0$ is fixed, we obtain solutions of
\GL\ and \GLhk\ as global minimizers for
$\EE$ and $\EEhk$ (in the appropriate spaces, $Y$ and $Z$): 
\begin{th}\label{existence}
For every fixed $g>0$, $\kappa\in \RR$, $d\in {\bf Z}-{0}$, the functional
$\II$  admits a minimizer $(u,v,m)\in X\times X \times H$.   Moreover, 
$(f,S,m)=(f_0+u,S_0+rv,m)$
is a smooth solution of the system \GL.
\end{th}

\begin{th}\label{existencehk}
For every fixed $g>0$ and $d\in {\bf Z}-{0}$, the functional
$\IIhk$  admits a minimizer $(u,m)\in X \times H$.  
Moreover, $(f,m)=(\finf+u,m)$
is a smooth solution of the system \GLhk.
\end{th}
The proofs of Theorems~\ref{existence}  and Theorem~\ref{existencehk}
are straightforward but technical, and are deferred to Section 6.

\subsection{Admissible solutions}

\medskip

As in \cite{ABG}, we define a natural class
of solutions to the  system \GL:
\begin{de}\label{admissible}
We call $(f_*, S_*, m_*)$ an \underbar{admissible} solution to
\GL\ if:
\begin{enumerate}
\item  \GL\ holds for all $r\in (0,\infty)$;
\item  $\EE (f_*,S_*,m_*)<\infty$;
\item  $f_*(r)\ge 0$ and $m_*(r)\ge 0$ for all $r\ge 0$;
\item  $S_*(0)=0$ and $m'_*(0)=0$.
\end{enumerate}
A solution $(f_*,m_*)$ of \GLhk\ is called \underbar{admissible} if
the above conditions hold, where we replace $\kappa$
by $\infty$ and disregard  $S_*$.  
\newline
A solution to \GL\ or \GLhk\ with $m_*\equiv 0$
is called a \underbar{normal core} solution.
\end{de}
The admissible solutions are those which are physically relevant in
the context of the vortex core problem described in the introduction.
We note that the normal core solutions are \underbar{unique} for
$\kappa^2\ge 2d^2$:  see \cite{ABG} for the case $2d^2\le \kappa^2<\infty$
and \cite{CEQ} for $\kappa=\infty$.

We now present some properties of admissible solutions.
In the following, we will assume that $\k \in \RR \bigcup \{ \infty \}$, 
with the understanding that $S_*=0$ when $\k =\infty$.
\begin{pr}\label{properties}
Let $(f_*, S_*, m_*)$ be any admissible solution of \GL.
Then:
\begin{enumerate}
\item  For all $r\in (0,\infty)$ it holds $0<f_*(r)<1$, $0\le m_*(r) < 1$, 
$f_*^2(r)+m_*^2(r)<1$,
and,  if $\k \neq \infty$, $0<S_*(r)<d$.
\item  Either $m_*(r)>0$ for all $r\in [0,\infty)$, or $m_*$ vanishes identically.
\item  $f_*(r)\to 1$, $m_*(r)\to 0$, and, if $\k \neq \infty$, $S_*(r)\to d$ 
as $r\to\infty$.  Moreover,
there exist constants $\sigma, C_0>0$ such that for $\k \neq \infty$
$$   
0<1-f_*(r)\le C_0 e^{-\sigma r}, \quad
0<d-S_*(r)\le C_0 e^{-\sigma r}, \quad 0\le m_*(r)\le C_0 e^{-\sigma r}, 
$$
and for $\k = \infty$
$$  
 0<1-f_*(r)\le {d^2\over {2 r^2}} + {{8d^2+d^4}\over{8r^4}} + {\cal O}(r^{-6}), \quad 0\le m_*(r)\le C_0 e^{-\sigma r}, 
$$
for all $r>0$.
\item  $f_*(r)\sim r^d$, $S_*(r)\sim r^2$ for $r\sim 0$.
\item  If $\k \neq \infty$, $S_*'(r) >0$ for all $r>0$.
\end{enumerate}\end{pr}
\Pf
The proof is very similar to that of Proposition~2.3 of \cite{ABG}, so we
provide only a sketch. 
From the  finiteness of the free energy we
immediately conclude that $m_* \in H$,
and hence $m_* \in L^p_r$, 
for any $p\in[2,\infty]$, and $m_*(r) \to 0$ as $r\to\infty$.
Since $f_*\ge 0$, finiteness of energy again implies $1-f_*\in L^2_r$ 
(see (\ref{a1}) for details,) and therefore the bound $0<f_*(r)<1$ follows
exactly as in Proposition~2.3 of \cite{ABG}.  When $\kappa<\infty$, the
bound $0<S_*(r)<d$ and the
proof that $S'_*(r)>0$ are also unchanged from \cite{ABG}.
To show $z=f_*^2+m_*^2 < 1$ we use the equation satified
by $z$:  this argument is already presented in \cite{ABBG}.
Statement (ii) is a simple consequence of the strong maximum
principle.

The exponential decay in (iii)  for $m_*$ is consequence of Proposition~7.4 in 
Jaffe \& Taubes \cite{JT}, and so are the ones for $f_*$ and $S_*$ if 
$\k \neq \infty$. If $\k =\infty$, the polynomial decay of $f_*$  can be proven 
as in Lemma 3.3 in \cite{CEQ},  since $m_*(r) \leq {C(R) \over {r^6}}$ for any $r>R$
with $C(R)$ a big enough constant.

The behavior at zero given in (iv) can be proven as 
in \cite{Plohr}. 
\QED

\medskip

We now connect admissible solutions to our space $X$.

\begin{pr}\label{connect}
Let $(f_0,S_0,m_0)$, $(f_1,S_1,m_1)$ be admissible solutions to
\GL.  Then $(f_1-f_0)\in X$, $[(S_1-S_0)/r]\in X$, and $m_1, m_0 \in H$.
\end{pr} 
\Pf As already remarked, condition (ii) of the 
definition of admissible solutions implies 
$m_1, m_0 \in H$, and $m_1, m_0 \in L^p_r$ for any $p \in [2,\infty)$. 
Then, the rest of the proposition for $\k \neq \infty$ is proven as in 
Proposition 2.4 of \cite{ABG}. When $\k=\infty$, we note that 
$(1-f_i) \in H$ for $i=1,2$, and that by (iv) of Proposition~2.5 
we  have $(f_1-f_2)^2 \leq c r^{2d}$ for $r\sim 0$ and again by finiteness of 
energy we conclude our statement.
\QED

\begin{re}\label{newdef}\rm
In light of Proposition~\ref{connect} we observe that
the choice of $f_0,S_0$ in the definition of the space
$Y$ may be replaced by any fixed admissible solution of the
equations \GL.  It will be convenient to choose instead the ``basepoint''
$(\tilde f_\kappa,\tilde S_\kappa,0)$ to be a ``normal core'' solution to
\GL.  In other words, an equivalent definition of the space
$Y$ is:
\be\label{Ydef}
Y=\{(f,S,m): \ f=\tilde f_\kappa + u, \ S=\tilde S_\kappa+rv, \
               u,v\in X, \  m\in H\}
\ee
We recall that the normal core solutions are uniquely determined
for $\kappa^2\ge 2d^2$.  When $\kappa^2<2d^2$ we fix any
one.
\end{re}

\begin{re}\label{solution}\rm
Proposition~\ref{connect} also implies that the admissible
solutions are exactly those which arise from minimization
problems for $\EE$ 
and $\EEhk$ in the space $Y$.  In particular, as an immediate
corollary we obtain the following statement:

\noindent
$(f_*,S_*,m_*)$  is an admissible solution to 
\GL\  if and
only if $f_*\ge 0$, $m_*\ge 0$, $(f_*,S_*,m_*)\in Y$ 
 and $\EE'(f_*,S_*,m_*)[u,v,w]=0$ 
 for all $u,v\in X$, and $w \in H$.

\noindent
An analogous statement holds for the problem \GLhk.
\end{re}

\medskip

With this choice of representation for our spaces $Y$, $Z$,
we now look at the second variation of energy with
respect to the variables $(u,v,w)\in X\times X\times H$.
We define
\bea\label{secvar}
\EE''(f_*,S_*,m_*)[u,v,w] &=&
  \left. {d^2\over dt^2}\right|_{t=0}
      \EE(f_*+tu, S_*+trv,m_*+tw) \nnn \\
&=&  \int\left\{
       (u')^2 +  (w')^2 +  {(d-S_*)^2\over r^2}u^2 +  \kappa^2 g w^2 +(v')^2  \right.\\
\nnn
&&\qquad  
  + {v^2\over r^2} -4 {(d-S_*)\over r} f_* uv  + f_*^2 v^2 \\
\nnn
&&\qquad \left.
 - \kappa^2 (1-f_*^2-m_*^2)(u^2 + w^2) + 2\kappa^2  (f_* u + m_* w)^2 \right\}\rdr. \\
\label{secvarinf}
\EEhk''(f_*,m_*)[u,w] &=&
 \left. {d^2\over dt^2}\right|_{t=0}
      \EEhk(f_*+tu,m_*+tw) \\
\nnn
&=&  \int\left\{
       (u')^2 +  (w')^2 +  {d^2\over r^2}u^2 +  g w^2 \right.\\
\nnn
&&\qquad \left.
 -  (1-f_*^2-m_*^2)(u^2 + w^2) + 2  (f_* u + m_* w)^2 \right\}\rdr.
\eea
Note that if we write $f_*=\fk + u_*$, $S_*=\Sk + rv_*$,
then 
$$\EE''(f_*,S_*,m_*)[u,v,w]= D^2\II(u_*,v_*,m_*)[u,v,w],  $$
the usual second Fr\'echet derivative.

\medskip

For admissible solutions which are stable, in the sense that
the second variation of energy about the solution is a non-negative
quadratic form, we have monotonicity of the profiles $f(r)$, $m(r)$.
\begin{th}\label{monotonicity}
Suppose $(f,S,m)$ is an admissible solution of \GL, and 
$\EE''(f,S,m)\ge 0$ as a quadratic form acting on $Y_0$.
Then $f'(r)>0$ and (if it is not identically
zero) $m'(r)<0$ for all $r>0$.
\end{th}
For the  problem \GLhk\ the same theorem holds, with exactly
the same proof.  We will see later that all admissible solutions 
of \GLhk\ with
$m(r)>0$ are stable (in the above sense), and hence we will obtain
the stronger result announced in Corollary~\ref{monoinf}.

\Pf
Let $\tilde u(r) = f'(r)$, $\tilde w(r)= m'(r)$.  Then, differentiating the
first and third equations of \GL,
\beann
&&  -\tilde u'' - {1\over r}\tilde u' + {(d-S)^2\over r^2} \tilde u
            -\kappa^2(1-3f^2-m^2)\tilde u + 2\kappa^2 mf\tilde w \\
&& \qquad \qquad
   = -{1\over r^2} \tilde u + 2 {d-S\over r} f \left[ {S'\over r} + {d-S\over r^2}\right] ,\\
&&  -\tilde w'' - {1\over r}\tilde w' + g\kappa^2 \tilde w
            -\kappa^2(1-f^2-3m^2)\tilde w + 2\kappa^2 f m \tilde u =
    -{1\over r^2} \tilde w.
\eeann
Suppose there exist intervals $(a,b)$, $(c,d)$ such that
\beann
&& \tilde u(r)<0 \quad r\in (a,b), \quad \tilde u(a)=0=\tilde u(b); \ \mbox{or} \\
&&  \tilde w(r)>0 \quad r\in (c,d), \quad \tilde w(c)=0=\tilde w(d).
\eeann
Note that by the properties (i), (iii) and (iv) of admissible
solutions in
Proposition~\ref{properties}, $a\neq 0$, $b,d <+\infty$.
Let 
$$
u(r)=\cases{  \tilde u(r), &  if $r\in (a,b)$,\cr
        0, &  otherwise,\cr  } \qquad
w(r)=\cases{  \tilde w(r), &  if $r\in (c,d)$,\cr
        0, &  otherwise,\cr  }.
$$
Then $u\le 0$, $w\ge 0$, and an integration by parts shows that
$$   \int (u')^2\rdr = -\int_a^b \tilde u \, {1\over r}(r\tilde u')'\rdr, $$
and similarly for $w$.  If we now use $(u,0,w)$ as a test function
in the second variation of energy and recall from Proposition~\ref{properties}
that $S(r)<d$, $S'(r)>0$ for all $r>0$, we obtain
\beann
0 &\le & \EE''(f,S,m)[u,0,w] \\
 & = & \int \left[ -{1\over r^2} u^2 + 
    2 {d-S\over r} f \left[ {S'\over r} + {d-S\over r^2}\right] u
         - {1\over r^2} w^2 \right] \rdr  \ 
< 0,
\eeann
unless $u,w\equiv 0$.  Consequently, $\tilde u = f'\ge 0$ and $\tilde
w=m' \le 0$.  Strict inequality follows from the Strong Maximum
Principle, since $\tilde u, \tilde w$ satisfy equations of
the form
\beann
&&  -\Delta_r \tilde u + c_1(r) \tilde u 
       \ge -2\kappa^2m f \tilde w \ge 0, \\
&&  -\Delta_r \tilde w + c_2(r) \tilde w 
       = -2\kappa^2m f \tilde u \le  0.
\eeann
\QED

\section{Nondegeneracy of solutions of \GLhk}

\setcounter{equation}{0}
\setcounter{th}{0}

\begin{th}\label{nondeg}
  For any admissible solution $(f_*,m_*)$ 
of \GLhk\ with $m_*>0$ there exists a constant $\sigma_*>0$ such that
$$      \EEhk'' (f_*,m_*)\, [u,w]  \ge \sigma_* (\|u\|_X^2 + \|w\|_H^2),  $$
for all $u\in X$, $w\in H$.
\end{th}
\begin{co}\label{monoinf}
For any admissible solution $(f_*,m_*)$ of \GLhk, $f'_*(r)>0$ for all $r\ge 0$.
If $m_*$ is not identically zero, then $m'_*(r)<0$ for all $r>0$.
\end{co}
The corollary follows from Theorem~\ref{nondeg} and the argument of
Theorem~\ref{monotonicity} when $m_*>0$.  Note that when $m_*\equiv 0$
the system \GLhk\ reduces to the single equation studied in
\cite{CEQ} and the strict monotonicity of $f_*$ is
part of their result. Also, in the case that $m_*\equiv 0$
the Theorem reduces to
$\EEhk'' (f_*)\, [u]  \ge \sigma_* \|u\|_X^2 $.

The key step in proving Theorem~\ref{nondeg} is the following
identity:

\begin{th}\label{tg}
For any admissible solution $(f_*,m_*)$
of \GLhk\  with $m_*>0$, and any $u\in X$, $w\in H$,
\be\label{lowerbd}
\EEhk'' (f_*,m_*)[u,w] =  \int
     \left\{   
        f_*^2 \left[ \left( {u\over f_*}\right)' \right]^2 
        +   m_*^2 \left[ \left( {w\over m_*}\right)' \right]^2 
      + 4 (f_* u + m_* w)^2
         \right\} \rdr
\ee
\end{th}
{\bf Proof of Theorem~\ref{tg}:} \
First we prove the identity for  $u\in C_0^\infty((0,\infty))$
and $w\in C_0^\infty([0,\infty))$.
First, note that using $f_*>0$ and $m_*>0$, we have
\be\label{identity}   f_*^2 \left[ \left( {u\over f_*}\right)' \right]^2 
    = (u')^2 - 2{u u' f'_*\over f_*} + u^2 { (f'_*)^2\over f_*^2} ,  
\ee
with a similar identity holding for $m_*, w$.
Hence,
\beann
0 & = &
\EEhk'(f_*,m_*) \left[  {u^2\over f_*},{w^2\over m_*}  \right] \\
& = &
\int \left\{
   (u')^2  + (w')^2 + {d^2\over r^2}u^2 + gw^2 - (1-f_*^2-m_*^2)(u^2+w^2) 
          \right.  \\
& & \qquad \left.
    -f_*^2 \left[ \left( {u\over f_*}\right)' \right]^2 
      - m_*^2 \left[ \left( {w\over m_*}\right)' \right]^2 \right\} \rdr.
\eeann

Substituting this in the formula 
 for $\EEhk''(f_*,m_*)[u,w]$
we obtain
$$
\EEhk'' (f_*,m_*)[u,w] = \int\left\{
     f_*^2 \left[ \left( {u\over f_*}\right)' \right]^2 +
       m_*^2 \left[ \left( {w\over m_*}\right)' \right]^2
    + 2 [f_*u + m_*w]^2 \right\} \rdr.
$$

To obtain the result for any $(u,w)\in X\times H$, let
$u_n$ be a sequence of $C_0^\infty((0,\infty))$ functions converging to
$u$ in $X$, and $w_n$ a sequence in
$C_0^\infty([0,\infty))$ converging to $w$ in $H$.
 By continuity of $\EEhk''(f_*,S_*)$, the limit
passes in the second variation of $\EEhk$. For the right hand side we expand,
$$ \label{quotient}
\int   f_*^2  \left(  \left(  {u\over f_*}  \right)'  \right)^2  \rdr
=\int \left\{  (u')^2-2{f'_*\over f_*} u u' + 
           \left({f'_*\over f_*}\right)^2   u^2   \right\}\rdr,
$$
and note that
$$\left({f'_*\over f_*}\right)^2\le c(1+{1\over r^2})$$
since $f_*\sim r^d$ for $r\sim 0$. Hence each term is controlled by
the $X$-norm and can be passed to the limit.  A similar argument may
be applied for the second term in the right-hand side of (\ref{lowerbd}).
The quotient is expanded as in (\ref{quotient}) above, with $m_*,w$
replacing $f_*,u$.   Then we claim
 that $m'(r)/m(r)$ is uniformly bounded for $r\in [0,\infty)$.
Indeed, by the basic gradient bound for solutions of
the Poisson equation (see section 3.4 of \cite{GT})
we have for any $r_0>1$,
$$
  |m'(r_0)| \le
    2\sup_{|r-r_0|\le 1} m(r)  + 
     \frac12 \sup_{|r-r_0|\le 1} | \kappa^2 (1-g-f^2-m^2)m | 
 \le C_1 \sup_{|r-r_0|\le 1} m(r).
$$
Applying the Harnack inequality (Corollary~9.25 of \cite{GT}) 
we then obtain:
$$
\left| {m'(r_0)\over m(r_0)}\right|
   \le   C_1 {\sup_{|r-r_0|\le 1} m(r) \over m(r_0)}
        \le C_1  {\sup_{|r-r_0|\le 1} m(r) \over  \inf_{|r-r_0|\le 1} m(r)}
            \le C'_1,
$$
for all $r_0>1$.  Therefore $m'/m$ is uniformly
bounded, and we may pass to the $H^1_r$ limit in  the second term in
(\ref{lowerbd}). The last term is clearly continuous in the $L^2_r$-norm in both $u$
and $w$. In conclusion, we may pass to the limit $u_n\to u$, $w_n\to w$ and
obtain (\ref{lowerbd}) for $u\in X$, $w\in H$.
\QED

\vskip 7pt

\noindent
{\bf Proof of Theorem~\ref{nondeg}:}\
Define
$$ \sigma_* = \inf\{
      \EEhk'' (f_*,m_*)[u,w]: \ 
              u\in X, \ w\in H, \ 
                 \|u\|_X^2 + \|w\|_H^2 =1\}.  $$
We must show that $\sigma_*>0$.

By Theorem~\ref{tg}, $\sigma_*\ge 0$.  To obtain
a contradiction, assume instead that $\sigma_*= 0$.  We claim that in this case 
 the infimum is attained at a nontrivial $(u_*,w_*)$, with
$\EE''(f_*,m_*)[u_*,w_*]=\sigma_*=0$.  But this contradicts Theorem~\ref{tg},
and hence $\sigma_*>0$.

We now claim that the infimum $\sigma_*=0$  is attained in $Z_0$.
Take any minimizing sequence: $(u_n,w_n)\in X\times H$ with
$\|u_n\|_X^2 + \|w_n\|_H^2 =1$ and 
$$   \EEhk'' (f_*,m_*)[u_n,w_n] \to \sigma_*=0.  $$
By the Sobolev embedding, there exists a subsequence
(still denoted by $u_n,w_n$) and $u_*\in X$, $w_*\in H$ so that
$u_n\to u_*$, $w_n\to w_*$, weakly in $X$, $H$ (respectively), and
 strongly in $L^2_{loc}$.

First, we claim that $(u_*, w_*)\neq (0,0)$.
Indeed, if both $u_*, w_*$ vanish identically then
by weak convergence $(u_n,w_n)\weak (u_*,w_*)=(0,0)$ and the
compact embeddings,
\beann
&&  \int \left(
    (u_n')^2 + {d^2\over r^2} u_n^2  + 2\kappa^2 f_*^2 u_n^2
        + (w'_n)^2 + g\kappa^2 w_n^2 \right)\, r\, dr \\
&& \qquad = 
         \EEhk'' (f_*, m_*)[u_n,w_n]\\
&& \qquad \qquad  + 
           \int \left[
              \kappa^2 (1-f_*^2-m_*^2)(u_n^2 + w_n^2) - 2\kappa^2 m_*^2 w_n^2
                 - 4\kappa^2 f_* m_* u_n w_n \right]\rdr \\
&&  \qquad \qquad \qquad
     \longrightarrow  0.
\eeann
In particular,
 $(u_n,w_n)\to (0,0)$
in the norm on $X\times  H$, which contradicts the fact that
$\|u_n\|_X^2 +\|w_n\|_H^2 =1.$  Thus the claim holds, and $(u_*,w_*)\neq (0,0)$.

Next, we use lower semicontinuity in the norm and
$L^2_{loc}$ convergence to pass to the limit,
\be \label{e1}
 \EEhk'' (f_*,m_*)[u_*,w_*] \le \liminf_{n\to\infty}
    \EEhk'' (f_*,m_*)[u_n,w_n] =0.
\ee
This contradicts 
Theorem~\ref{tg}, since $\EEhk'' (f_*,m_*)[u_*,w_*]>0$.
(Note that $u/f_*$ is non-constant since $u\in X$ but $f_*\not\in X$.)
We conclude that $\sigma_*>0$, as desired.
\QED

We note that the same result holds when $m_*\equiv 0$. Hence following the method
of \cite{ABG}, we obtain another proof of uniqueness for the solution to the 
high kappa equation for $f_*$ studied in \cite{CEQ}.

\section{Bifurcation from the normal cores}

\setcounter{equation}{0}
\setcounter{th}{0}

In this section we show that (when $\kappa^2\ge 2d^2$)
AF core solutions are nucleated
by means of a bifurcation from the normal core solution family
at a simple eigenvalue of the linearized equations.   We will also
require {\it a priori} estimates (whose proof we will present
in Section~6)  to obtain
global information about the solutions set for all $\kappa^2\ge 2d^2$,
and the stronger result of Theorem~\ref{nondeg} to fully categorize solutions  
in the extreme type-II model \GLhk.
We present the detailed argument for the problem \GL.  The functional analytic
framework is entirely similar for the problem \GLhk, and so we
omit it and concentrate instead on the more precise global charaterization
of solutions which we prove for \GLhk.

\subsection{Local bifurcation at $g^*_\kappa$}

We define a map $\FF: \  Y\times \RR\to Y_0^*$
by
$$  
\langle (u,v,w), \FF (f_*,S_*,m_*,g)\rangle_{Y_0,Y_0^*}
     = \EE' (f_*,S_*,m_*)[u,v,w], 
$$
$(u,v,w)\in Y_0$, $(f_*,S_*,m_*)\in Y$.  Its linearization
is the operator $\FF'(f_*,S_*,m_*,g)\in L(Y_0,Y_0^*)$ defined
by
$$  \langle (u,v,w), \FF'(f_*,S_*,m_*,g)[\phi,\psi,\xi]\rangle_{Y_0,Y_0^*}
    = \left. {d\over dt} \right|_{t=0}
          \EE'(f_*+t\phi,S_*+rt\psi,m_*+t\xi) [u,v,w].  $$
We remark that  the explicit expansion of the
energy $\II$ in terms of $u_*=f_*-\fk$, $v_*=(S_*-\Sk)/r$, $w_*$
ensures that $\FF$ is a $C^2$ map in all arguments $u_*,v_*,w_*,g$.

By the natural identification $Y_0\simeq Y_0^*$ of
a Hilbert Space with its dual, we may also represent
$\FF'$ by $\LL_g\in L(Y_0,Y_0)$ as
$$  \left( (u,v,w), \LL_g [\phi,\psi,\xi] \right)_{Y_0}
         =  \langle (u,v,w), \FF'(f_*,S_*,m_*,g)[\phi,\psi,\xi]\rangle_{Y_0,Y_0^*}.
$$
If $i: Z^*\to Z$ is the isomorphism, then $\LL_g= i \circ \FF'(f_*,S_*,m_*,g)$.

\begin{lem}\label{Fredholm}
For all $g>0$, $\LL_g$ is a Fredholm operator of index zero.
\end{lem}
\Pf
Define an equivalent inner product on $Y_0$,
\beann
\left( (u,v,w), (\phi,\psi,\xi)\right)_{Y_0}
  &= & \int \left\{  u'\phi' + 2\kappa^2 u\phi +
                 {(d-S_*)^2\over r^2} u\phi + v'\psi'  \right. \\
    && \quad  \left.
         + v\psi + {1\over r^2} v\psi + w'\xi' + g\kappa^2 w\xi\right\}
                         \rdr.
\eeann
Then we write
$$  \left( (u,v,w), \LL_g [\phi,\psi,\xi] \right)_{Y_0}
     = \left( (u,v,w), (\phi,\psi,\xi)\right)_{Y_0}
        +     \left( (u,v,w), K [\phi,\psi,\xi] \right)_{Y_0} ,
$$
where $K$ is defined by 
\beann
 \left( (u,v,w), K [\phi,\psi,\xi] \right)_{Y_0} &=&
\int \left[ 2\kappa^2(f_*^2-1)u\phi + 2\kappa^2 f_* m_* (u\xi + w\phi)
+ 2\kappa^2 m_*^2 w\xi   \right.\\
&& \qquad 
-\kappa^2(1-f_*^2-m_*^2)(u\phi + w\xi) + (f_*^2-1) v\psi \\
&&\qquad \left.
     - 2{d-S_*\over r} f_* (u\psi + v\phi) \right]\rdr.
\eeann
Recalling the decay properties of
$f_*,S_*,m_*$ and the embedding properties of $H,X$ we
observe that  $K$ is compact, and hence $\LL_g= Id_{Y_0} + K$
is Fredholm with index zero.
\QED

As a direct consequence of Lemma~\ref{Fredholm},
$$  \dim \ker(\FF')=\dim \ker(\LL_g)=\codim \Ran(\LL_g) = \codim\Ran (\FF').  $$

Now we may apply the standard bifurcation theory of
Crandall \& Rabinowitz \cite{CR} at an eigenvalue $g^*$  of $\FF'(\fk,\Sk,0, g^*)$. 
Indeed, note that when $m_*=0$ the linearization of $\FF$ decouples
into two components,
$$  \langle (u,v,w), \FF' (f_*,S_*,0,g)[\phi,\psi,\xi] \rangle_{Y_0,Y_0^*}
    =\langle (u,v), \FF'_{1,2}(f_*,S_*)[\phi,\psi]\rangle_{X^2,(X^2)^*}
          + \langle w, \FF'_3 (f_*,g) \xi \rangle_{H,H^*},  $$
where 
\beann
  \langle (u,v), \FF'_{1,2}(f_*,S_*)[\phi,\psi]\rangle_{X^2,(X^2)^*}
  &= &
\langle (u,v,0), \FF' (f_*,S_*,0,g)[\phi,\psi,0] \rangle_{Y_0,Y_0^*} \\
& = &  \int \left[
   u' \phi' + {(d-S_*)^2\over r^2} u\phi + v'\psi'
    + {v\psi\over r^2}  \right. \\
&& \ \left.
+ f_*^2 u\psi - 2{d-S_*\over r} f_* (u\psi + v\phi)
-\kappa^2 (1-3f_*^2) u\phi \right]\rdr,
\eeann
and
\beann
  \langle w, \FF'_3(f_*,g)\xi \rangle_{H,H^*}
  &= &
\langle (0,0,w), \FF' (f_*,S_*,0,g)[0,0,\xi] \rangle_{Y_0,Y_0^*} \\
& = &  \int \left\{ w'\xi' + g\kappa^2 w\xi -\kappa^2 (1-f_*^2)w\xi
   \right\}\rdr.
\eeann
By Theorem~3.1 of \cite{ABG}, when $\kappa^2\ge 2d^2$
the operator $\FF'_{1,2}\ge \sigma_*>0$ is
bounded away from zero (in quadratic form sense.)
Hence, if $(\phi,\psi,\xi)\in \ker (\FF'(\fk,\Sk,0,g^*_\kappa))$, we take 
$(u,v,w)=(\phi,\psi,0)$ and obtain
\beann
0 &=&  \langle (\phi,\psi,0), \FF' (f_*,S_*,0,g^*_\kappa)[\phi,\psi,\xi] \rangle_{Y_0,Y_0^*}\\
&=& \langle (\phi,\psi), \FF'_{1,2}(f_*,S_*)[\phi,\psi] \rangle_{X^2,(X^2)^*} \\
&\ge & \sigma_* (\|\phi\|_X^2 + \|\psi\|_X^2)
\eeann
In particular, $\phi,\psi=0$.  

The operator $\FF'_3(f_*,g)= L+g\kappa^2$
where $L =-\Delta_r - V(r)$  is a Schr\"odinger operator
with  potential
$V(r)=\kappa^2 (1-f_*^2(r)) \ge 0$ and $V(r)\to 0$ as $r\to\infty$.
It is a well-known fact in mathematical physics that in dimension two, such operators
have at least one negative eigenvalue:
\begin{lem}\label{boundstate}
Suppose $V: \ [0,\infty)\to\RR$ is continuous, non-negative,
$V(r)\to 0$ as $r\to\infty$, and $V$ is not identically zero,
and define $L=-\Delta - V(r)$ as a  self-adjoint operator on the
space $L^2(\RR^2)$.  Then the ground state energy,
$$  \lambda_0 = \inf \left\{
            { \int [(u')^2 - V(r) u^2] \rdr \over \int u^2\rdr} :\
                   u\neq 0, \  u\in H\right\} <0,  $$
and is attained at an eigenfunction $u_0\in H$.  Moreover, 
$\lambda_0$ is an isolated, non-degenerate eigenvalue,
$u_0\in H$, and $u_0>0$.
\end{lem}
The proof follows as an application of the Birman--Schwinger principle
in Reed \& Simon \cite{RS}.  We provide an elementary variational
proof for the reader's convenience.

\Pf
Let 
$$  u_n(r)= \cases{
       1, &  if $r\le n$, \cr
      {ln(r/n^2)\over ln(1/n)}, & if $n\le r\le n^2$, \cr
       0, & if $r\ge n$.\cr  }
$$
Then,
$$  \int (u'_n)^2 \rdr = {1\over \ln n} \to 0,  $$
while
$$  \int V(r) u_n^2\rdr  \ge \int_0^n V(r)\rdr \to \int_0^\infty V(r)\rdr >0 $$
(possibly infinite.)
Hence, for $n=N$ large but fixed we have $\int [ (u'_N)^2 - V(r) u_N^2 ]\rdr <0$,
and hence $\lambda_0<0$.
Since $L$ is a relatively compact perturbation of $-\Delta$, $\lambda_0$
is a discrete eigenvalue with associated eigenfunction $u_0$ contained in the form
domain of $L$, $H$.  By standard arguments, $u_0>0$ and 
$\lambda_0$ is a simple (non-degenerate) eigenvalue.
\QED

By Lemma~\ref{boundstate} 
$$  -g_\kappa^* = \inf_{w\in H-\{0\}}  
    {\int \left[ {1\over \kappa^2} (w')^2 - (1-\fk^2) w^2\right]\rdr
        \over   \int w^2\rdr } < 0,
$$
and $\lambda_0=-\kappa^2 g_\kappa^*$
is the ground state eigenvalue of the Schr\"odinger operator
$-\Delta_r -\kappa^2(1-\fk^2)$.   Since $\lambda_0$ is a simple eigenvalue
$$  \dim\ker (\FF'_3 (\fk,g^*_\kappa)) = 
    \dim\ker (-\Delta_r - \kappa^2(1-\fk^2) +\kappa^2 g_\kappa^*) = 1.  $$
In conclusion   when $g=g_\kappa^*$ the operator $\FF'(\fk,\Sk,0,g^*_\kappa)$ has
a simple eigenvalue and  the eigenvector is of the form $(0,0,w_\kappa)$
with $w_\kappa$ the (positive) eigenfunction of $\FF'_3$.

Finally, we observe that the operator 
$({\partial\over\partial g})\FF'(\fk,\Sk,0,g)\in L(Y_0,Y_0^*)$,
$$  \langle (u,v,w), {\partial\over\partial g}\FF' (\fk,\Sk,0,g) [\phi,\psi,\xi]\rangle_{Y_0,Y_0^*}
     =  \int \kappa^2 w\xi\rdr.  $$
At the eigenvalue $g=g_\kappa^*$ we have
$$  {\partial\over\partial g}\FF' (\fk,\Sk,0,g^*_\kappa) [0,0,w_\kappa] = \kappa^2w_\kappa
   \not\in \Ran (\FF'(\fk,\Sk,0,g_\kappa)).  $$
Therefore Theorem~1.7 of \cite{CR} applies, and $g_\kappa^*$ is
a bifurcation point for $\FF$ in $Y\times \RR$:
there exists a neighborhood $U$ of $(\fk,\Sk,0,g^*_\kappa)$ in $Y_0\times \RR$,
such that the set of non-trivial solutions of $\FF (f,S,m,g)=0$ in $U$ is 
a unique $C^1$ curve parametrized by $\ker (\FF' (\fk,\Sk,0, g^*_\kappa))$. 

\begin{re}\label{bifdir}\rm
Since $\FF$ is a smooth ($C^\infty$) map, we may calculate various
derivatives of the bifurcation curve through the 
normal core solutions at $g^*_\kappa$.  If we parametrize $g=\gamma(t)$,
with $\gamma(0)=g^*_\kappa$, then we follow 
Crandall \& Rabinowitz \cite{CR} or Ambrosetti \& Prodi \cite{AP}
(see Remarks 4.3) to calculate derivatives of $\gamma(t)$ and determine
the direction of the bifurcation curve locally at $g=g_\kappa^*$.
We obtain that $\gamma'(0)=0$, and 
$$  \gamma''(0)= -2{
        \int \left[ \fk u_* w_\kappa^2 +w_\kappa^4 \right]\rdr
               \over
      \int w_\kappa^2 \rdr },
$$
where $u_\kappa$ is obtained from the (unique) solution to the linear system
$$
\FF'(\fk,\Sk,0,g_\kappa^*) [u_*, v_*,w_*] 
     = -\left( 2\kappa^2 \fk w_\kappa^2, 0, 0 \right)  $$
with $(u_*,v_*,w_*)\perp \ker \FF'(\fk,\Sk,0,g_\kappa^*)$.  
By taking the scalar product of the above system with $(u_*,v_*,0)$
(and recalling that $\FF'(\fk,\Sk,0,g_\kappa^*)$ is positive definite
in the complement of its kernel)
we obtain $\int \fk u_* w_\kappa^2\rdr < 0$, and hence the
expression for $\gamma''(0)$ is indefinite in sign.  In a joint
paper with J. Berlinsky \cite{ABBG} we present computational
evidence that solutions bifurcate to the left, to smaller values
$g<g_\kappa^*$.   By standard bifurcation theory (see \cite{CR2},
for example) the direction of bifurcation indicates the stability of the
solutions, and indeed we observe numerically that the AF core
solutions which bifurcate at $g_\kappa^*$ are stable (local energy
minimizers.)  
\end{re}

\subsection{Global bifurcation for \GLhk}

We obtain the same abstract bifurcation result for the extreme type-II model, \GLhk.
Namely, the value 
$$  g_\infty^* = -\inf_{w\in H-\{0\}}  
    {\int \left[  (w')^2 - (1-\finf^2) w^2\right]\rdr
        \over   \int w^2\rdr }> 0
$$
is a bifurcation point for nontrivial ($m>0$) solutions from the
(trivial) curve of normal core solutions $(\finf,0,g)$.
But in this case we can make a much more precise statement:

\begin{pr}\label{graph}
Let
$$  \Sigma = \{ (f,m,g) : \ 
           \mbox{ $(f,m)$ is an admissible solution to \GLhk\ with $m>0$}\}.  $$
Then ${\cal C}=\Sigma\cup \{(\finf,0,g_\infty^*)\}$ is a connected $C^1$
curve, parametrized
by $g$.  Moreover for any $g_0>0$, ${\cal C}\cap \{ g\ge g_0\}$ is compact.
\end{pr}
As a consequence we have the following exact solvability theorem
for \GLhk.
\begin{th}\label{unique}
For $g\ge g_\infty^*$, the normal core solutions $(\finf,0)$ are
the only admissible solutions of \GLhk.
\newline
For $0<g<g_\infty^*$ there is a unique solution with $m>0$.  This solution
is the global minimizer of $\EEhk$. 
\end{th}
The proofs of these two results hinge on the powerful Theorem~\ref{nondeg}
and the following compactness theorem, which will be proven in Section~6:
\begin{th}\label{compactness}
Let $0<a<b$.  Then the set of all admissible solutions of \GLhk\
with $g\in [a,b]$ is compact in $Z$.
\end{th}

\noindent
{\bf Proof of Theorem~\ref{graph}:}\
Let $\cal C'$ be a maximally connected component of $\cal C$, and
suppose $(f_0,m_0,g_0)\in {\cal C'}$ but $(f_0,m_0,g_0)\neq (\finf,0,g_\infty^*)$.
Since $m=0$ only when $g=g_\infty^*$, we must have $m_0>0$.
By  Theorem~\ref{nondeg}, $(f_0,m_0,g_0)$ is a nondegenerate zero
of $\FF$ in $Z\times \RR$, so by the Implicit Function Theorem there exists
a neighborhood $U$ of $(f_0,m_0,g_0)$ in $Z\times \RR$, an
interval $J=(g_0-\delta,g_0+\delta)$, and a $C^1$
function $\Phi: J\to Z$ so that  all solutions
of $\FF=0$ in $U$ are of the form $(\Phi(g),g)$ with $g\in J$.

Let 
$$
  \hat g = \sup \{g:\  \mbox{ there exists a solution  $(f,m,g)\in {\cal C'}$} \} > g_0.
$$
Note first that any solution must satisfy
$$  0\le \int  (m')^2 \rdr \le \int (1-g)m^2 \rdr,  $$
and hence  $g<1$ for any solution with $m\not\equiv 0$.
Since by Proposition~\ref{compactness},
${\cal C}'\cap\{g\ge g_0\}={\cal C}'\cap\{g_0\le g\le 1\}$ is compact, 
there exists a solution at $g=\hat g$,
$(\hat f, \hat m ,\hat g)\in {\cal C'}$.   First, we claim that
$\hat m=0$.  If not, then by Proposition~\ref{properties} $\hat m(r)>0$ for all $r>0$,  
so by Theorem~\ref{nondeg}, $(\hat f,\hat m)$ is a nondegenerate minimum
of (GL)${}_{\hat g,\infty}$.  By the Implicit Function Theorem argument
above there exist a $C^1$ curve of nontrivial solutions through
$(\hat f, \hat m ,\hat g)$, parametrized by $g$.  In particular, we
contradict the definition of $\hat g$ is the supremum of all $g$ for
solutions in the connected component $\cal C'$.  Hence $\hat m=0$, as desired.

Now we show that $\hat g= g^*_\infty$.  Take a sequence $(f_n,m_n,g_n)\in {\cal C'}$
with $g_n\to \hat g$, so the above arguments imply that
$f_n-\finf\to 0$ in $X$   and $m_n\to 0$ in $H$.
Let
$$  t_n = \int (1-f_n^2) m_n^2 \rdr \to 0.  $$
Then $w_n=m_n/t_n$ solves
\be\label{weqn}  -w''_n - {1\over r} w'_n + g_n w_n = (1-f_n^2 -m_n^2)w_n.  
\ee
Since 
\be  \int \left(  (w'_n)^2 + g w_n^2 \right)\rdr
    = \int (1-f_n^2-m_n^2) w_n^2\rdr \le \int (1-f_n^2)w_n^2\rdr =1, 
\ee
(by the choice of $t_n$,) we have $\|w_n\|_H\le 1/g$ and we may
extract a subsequence (which we continue to call $w_n$) which
converges $w_n\weak w_\infty$ weakly in $H$ and strongly in
$L^2_{loc}$.  By the strong convergence of $f_n\to \finf$ we have
$\int (1-\finf^2)w_\infty^2 \rdr =1$, so $w_\infty\not\equiv 0$,
and $w_\infty \ge 0$.  Passing to the limit in (\ref{weqn})
we have 
\be\label{wweqn}
  \int \left(  w'_\infty \phi' +\hat g w_\infty\phi - (1-\finf^2)w_\infty\phi \right)\rdr
       =0,  
\ee
for all $\phi\in H$.  This can only occur when $\hat g=g^*_\infty$, the
ground state eigenvalue of the above Schr\"odinger operator.

We have just shown that the point $(\finf,0,g_\infty^*)$ belongs to every connected
component of $\cal C$, and hence $\cal C$ is connected.  
The solution set $\cal C$ is everywhere a $C^1$ curve:  for $g>g^*_\infty$
this results from the Implicit Function Theorem argument in the first
paragraph, and at $g_\infty^*$ it is a consequence of bifurcation from
a simple eigenvalue \cite{CR}.  
We now claim that there exists exactly one solution in $\cal C$ for
every $g\le g_\infty^*$.  Suppose not, and consider
\beann 
  &&  D= \{ g\in (0,g_\infty^*): \ 
            \mbox{there exist two distinct solutions $(f_{g,1},m_{g,1})$, 
             $(f_{g,2},m_{g,2})$
              in $\cal C$ at $g$.}\},  \\
&&\mbox{and } g_0 = \sup D.
\eeann
First, we note that $g_0<g_\infty^*$.  To see this we note that
the only solution in $\cal C$ with $g=g_\infty^*$ is the
normal core solution, and the bifurcation theorem 
ensures that the solution set in a neighborhood of the bifurcation
point $(\finf,0,g_\infty^*)$ is a single smooth curve.
 
Next, we claim that $g_0\not\in D$.  Indeed, if $g_0\in D$
there exist two distinct solutions $(f_{g_0,1},m_{g_0,1})$ and
$(f_{g_0,2},m_{g_0,2})$ for $g=g_0$.  By the Implicit Function Theorem
argument of the first paragraph there exist neighborhoods $U_1$ (of
$(f_{g_0,1},m_{g_0,1},g_0)$) and $U_2$ (of $(f_{g_0,2},m_{g_0,2},g_0)$)
in $Z\times\RR$ such that all solutions of ${\cal F}=0$ in 
$U_1$, $U_2$ are given by smooth curves parametrized by
$g$.  In particular, $\cal C$ contains two distinct solutions for
$g$ in an interval to the right of $g_0$, contradicting the definition
of $g_0$ as the supremum.

Hence $g_\infty^*>g_0\not\in D$, and there exists a sequence $g_k\to g_0$
for which ${\cal C}$ contains two distinct solutions, $(f_{g_k,1},m_{g_k,1})$,
$(f_{g_k,2},m_{g_k,2})$.  By Theorem`\ref{compactness}, along some subsequence
 these solutions converge, and since $g_0\not\in D$, they both converge
to a single solution, $(f_{g_0},m_{g_0})$.  But this contradicts the
Implicit Function Theorem argument, which implies that the solution
set near $(f_{g_0},m_{g_0}, g_0)$ is a single curve parametrized by $g$.
We conclude that the AF core solutions are {\it unique} for
each $g\in (0,g_\infty^*)$.
\QED

\subsection{Behavior for $g \to 0$, $\kappa<\infty$}

For the problem \GL\ we do not have the strong information provided
by Theorem~\ref{nondeg} which determines the global structure of the
solution set, and hence we cannot make the same elegant conclusion about
the uniqueness of AF core solutions.  However we
may still say something about the global structure of the continuum
bifurcating from the normal cores at $g=g_\kappa^*$.  When $\kappa^2\ge 2d^2$
we may apply the Global Bifurcation Theorem of Rabinowitz \cite{Rab2}
to conclude that the continuum $\Sigma_\kappa$ of zeros of
$\FF(f,S,m,g)=0$ with $m>0$ is unbounded in the space $Y\times \RR$.
(Note that $\Sigma_\kappa$ cannot contain any other eigenvalues of the linearization
about the normal core solutions, as is easily seen from the calculations 
(\ref{weqn})--(\ref{wweqn}) above.)  In the next section we will prove the
following {\it a priori}\/ estimate, which has as a direct consequence
the fact that $\Sigma_\kappa$ can only become unbounded as $g\to 0^+$:
\begin{th}\label{apriori}
Let $d,\kappa$ be fixed.  For any compact
interval $J\in (0,\infty)$ there exists $C_0=C_0(\kappa,d,J)>0$ such
that every admissible solution $(f,S,m)$ of \GL\ with $g\in J$
satisfies $\|(f,S,m)\|_Y\le C$.
\end{th}

Let us now concentrate on this loss of compactness in the continuum
$\Sigma_\kappa$ as $g\to 0^+$.  We prove:
\begin{th}\label{gtozero}
For any sequence of (absolute) minimizers  $(f_g,S_g,m_g) \in Y$ 
with $g\to 0^+$ we have  
$f_g\to 0$ in $X_{loc}$,  $S_g\to 0$ locally uniformly, and
$m_g\to 1$ in $H_{loc}$.
\end{th}
Fix $\k \in \RR$, and for any $g>0$ consider a minimizer $(f_g,S_g,m_g) \in Y$ of $\EE$.
\begin{lem}\label{Etozero}
$$
\EE (f_g,S_g,m_g) \longrightarrow 0, \qquad \mbox{ as } \qquad g \to 0.
$$
\end{lem}
\Pf
We will show that for any $\epsilon >0$ 
there exist $g_\epsilon > 0$ and $H$ 
radial functions $(f_\ep,S_\ep,m_\ep) \in Y$ such that 
$0<\EE(f_\ep,S_\ep,m_\ep) < \ep$, for any $g < g_\ep$.

For a fixed $\rho >0$, we define
$$  u_\rho(r)= \cases{
       1, &  if $r\le \rho$, \cr
      {\ln(r/\rho^2)\over \ln(1/\rho)}, & if $\rho\le r\le \rho^2$, \cr
       0, & if $r\ge \rho^2$,\cr  }
$$
we consider
$$
f_\rho(r)=\cos (u_\rho(r) \frac \pi 2), \qquad m_\rho(r)=\sin (u_\rho(r) \frac \pi 2),
$$
and
$$
\quad S_\rho(r)=\cases{  0, &  if $r \in (0,\frac \rho 2)$,\cr \cr
        d, &  if $r \in (\rho, \infty)$,\cr  }, 
$$
A direct computation shows that 
$$\EE(f_\rho,S_\rho, m_\rho) \leq \frac C{\rho^2}+
     \frac {\pi^2}{ 4\ln \rho} + \frac {\k^2} 2 g \rho^4$$ 
for any $g>0$. 

For a given $\ep > 0$, we choose a $\rhe$  such that 
$\frac C{\rho_\ep^2} + \frac {\pi^2} { 4\ln \rho_\ep}< \frac {\ep}2$, 
and a $g_\ep=g_\ep(\rhe)$ for which
$\frac {\k^2} 2 g_\ep \rho_\ep^4 < \frac {\ep}2$, i.e. $ g_\ep < \frac {\ep}2 
\frac 2{\k^2 \rho_\ep^4}$. Then,  
$\EE(f_g,S_g, m_g) \leq \EE(f_{\rhe},S_{\rhe}, m_{\rhe}) < \ep$, for any $g < g_\ep$.
\QED

{\bf Proof of Theorem~\ref{gtozero}}\
By Lemma~\ref{Etozero} each term in the energy tends to zero as $g \to 0$.
First, note that $\int (S'_g/r)^2\rdr \to 0$ combined with  (1.4) in \cite{BC}
implies that
\be\label{Sto0}
\mbox{$S_g(r)/r \to 0$ uniformly. }
\ee
For any $R_0>0$, we then have
\beann
o(1) &= & \int {(d-S_g)^2\over r^2} f^2_g \rdr
   \ge \int_0^{R_0} {(d-S_g)^2\over r^2} f^2_g \rdr  \\
&= & \int_0^{R_0} {d^2\over r^2} f^2_g \rdr  + o(1).
\eeann
In particular, $f_g\to 0$ in $L^2_{loc}$, $X_{loc}$.
Finally, by the reverse triangle inequality,
\beann
o(1) &=& \sqrt{ {\kappa^2\over 2} \int_0^{R_0} (1-f_g^2-m_g^2)^2\rdr }
    = {\kappa\over \sqrt{2}}\| 1-f_g^2-m_g^2\|_{L^2([0,R_0])} \\
& \ge & {\kappa\over \sqrt{2}} \left[\|1-m_g^2\|_{L^2([0,R_0])}
    - \|f_g^2\|_{L^2([0,R_0])} \right] \\
& \ge & {\kappa\over \sqrt{2}} \| 1-m_g\|_{L^2([0,R_0])} + o(1),
\eeann
where we have also used $0\le f_g< 1$, $0<m_g<1$, and $f_g\to 0$
in $L^2_{loc}$.  In conclusion $m_g\to 1$ in $L^2_{loc}$ and in fact
in $H^1_{loc}$, since $\int (m'_g)^2\rdr\to 0$ by the energy estimate.
\QED

\section{The limit $\kappa\to\infty$}

\setcounter{equation}{0}
\setcounter{th}{0}

In this section we show that the problem \GLhk\ arises as a
limiting case of \GL\ as $\kappa\to \infty$. 
For any solution $(f_\kappa,S_\kappa,m_\kappa)$ of
\GL, define
\be\label{rescale}
\fhk(r) = f_\kappa \left( {r\over \kappa}\right), \qquad
\Shk(r) = S_\kappa \left( {r\over \kappa}\right), \qquad
 \mhk(r) = m_\kappa \left( {r\over \kappa}\right).
\ee
We prove:
\begin{th}\label{limit}
Let $(f_\kappa,S_\kappa,m_\kappa)$ 
be any family of
solutions of \GL\ for $\kappa>0$,
and $(\fhk,\Shk,\mhk)$ defined as in (\ref{rescale}).  
For any sequence  $\kappa_n\to\infty$, there exists a subsequence
and a  solution
$(f_\infty,m_\infty)$ of \GLhk\ so that (as $\kappa_{n_k}\to\infty$,)
$\hat f_{\kappa_n}-f_\infty\to 0$
in $X$, $\hat m_{\kappa_n}-m_\infty\to 0$ in $H$, and 
$\hat S_{\kappa_n}\to 0$
locally uniformly.  Moreover:
\begin{enumerate}
\item  If $g\ge g_\infty^*$, then $m_{\kappa}\to 0$;
\item  If $m_\kappa\not\equiv 0$ for all large $\kappa$
and $g\neq g_\infty^*$,
then $\lim_{\kappa\to\infty} \mhk = m_\infty>0$.
\end{enumerate}
\end{th}
As a simple consequence of the uniform convergence
of $\fhk\to\finf$ we have the following
\begin{co}\label{glimit}
$$  g_\infty^* = \lim_{\kappa\to\infty} g_\kappa^* . $$
\end{co}
\begin{re}\label{limitremark}\rm
This implies that the bifurcation diagram for \GL\
with $\kappa$ very large should strongly resemble the
very precise image given for \GLhk\ by Theorem~\ref{unique}.
In particular, for any fixed $g>g_\infty^*$ 
\GL\ cannot have solutions $(f_{\kappa,g},S_{\kappa,g},m_{\kappa,g})$
with $m_{\kappa,g}>0$ for $\kappa$ large.
\end{re}

\medskip

Simple calculations using the energy $\EE$ show that 
$\inf_Y \EE \sim \ln\kappa$, and hence we require
 require energy-independent estimates for our
solutions $(\fhk,\Shk,\mhk)$.  To obtain these estimates we
begin with a simple version of the celebrated Pohozaev
identity.  This identity will also be essential for proving the
{\it a priori}\/ estimates used in the bifurcation analysis in
the previous section.
\begin{pr}\label{pohozaev}
For any finite energy solution $(f,S,m)$ of \GL\ we have
$$   g\kappa^2 \int m^2\rdr  
         + {\kappa^2\over 2}\int (1-f^2-m^2)^2\rdr
= \int \left[ {S'\over r} \right]^2 \rdr. $$
For any finite energy solution $(f,m)$ of \GLhk\ we have
$$    g \int m^2\rdr  
         + {1\over 2}\int (1-f^2-m^2)^2\rdr
= {d^2\over 2}.  $$
\end{pr}
\Pf
We multiply the first equation in \GL\ by $f'(r) r$ and integrate
$r \, dr$ to obtain:
\beann
  {\kappa^2\over 2}\int (1-f^2-m^2)(f^2)' r^2\, dr
      & = &
\int  (d-S)^2 \left( \frac12 f^2\right)' \, dr \\
& = & \int (d-S) S' f^2 \, dr =  -\int \left({S'\over r}\right)' S' \rdr\\
& = & \int {S'\over r} (S' r)'\, dr  = \int \left( {S'\over r}\right)^2 \rdr,
\eeann
using the equation for $S(r)$, and integrating by parts whenever 
necessary.  We also multiply the third equation in \GL\ by
$m'(r) r$ and integrate $r \, dr$ to obtain:
$$    {\kappa^2\over 2} \int (1-f^2-m^2)(m^2)' r^2\, dr
   = {g\kappa^2\over 2} \int (m^2)' r^2\, dr
         = -g\kappa^2 \int m^2 \rdr.  $$
Together,
\beann
\int \left( {S'\over r}\right)^2 \rdr
& = &  g\kappa^2 \int m^2 \rdr 
    + {\kappa^2\over 2} \int (1-f^2-m^2)(m^2+f^2)' \, r^2\, dr \\
& = & g\kappa^2 \int m^2 \rdr 
    + {\kappa^2\over 2} \int (1-f^2-m^2)^2\rdr.
\eeann

For the case $\kappa=\infty$ we proceed in the same way, except
the equation for $f$ yields
$$  \int (1-f^2-m^2) \left( {f^2\over 2}\right)' r^2\, dr
        = {d^2\over 2} .  $$
The calculation then continues as above.
\QED

\medskip

\Pf  {\bf Step 1:} Bounding the sequence.

\noindent
From the Pohozaev identity (Proposition~\ref{pohozaev}) and
Lemma~4.2 of \cite{BC} after rescaling we have:
\bea
\nnn
d^2 &\ge & \int \left( {S'_\kappa\over r}\right)^2 \rdr \\
\label{L1}
& = & \kappa^2 \int \left( {\Shk'\over r}\right)^2 \rdr \\
\label{L2}
&= & g\int \mhk^2 \rdr + {1\over 2} \int (1-\fhk^2 - \mhk^2)^2 \rdr.
\eea
Using (\ref{L1}) and Lemma~1.2 (ii) in \cite{BC} we have
\be\label{L5}
\sup_{r\in [0,\infty)} \left|{\Shk(r)\over r}\right| \to 0,
\ee
and hence
$\Shk \to 0$ locally uniformly.
From (\ref{L2}) we obtain the uniform bound
$\|\mhk\|_2\le C$ (depending on $g$, which we assume
is fixed.)
From the equation for $m_\kappa$, after a change of scale,
we obtain:
$$  \int  \left[  (\mhk')^2 + g \mhk^2 \right]\rdr = 
      \int (1-\fhk^2-\mhk^2)\rdr \le \int \mhk^2\rdr \le C, $$
and therefore $\|\mhk\|_H\le C$ uniformly in $\kappa$.

Recalling Proposition~\ref{properties}, any
solution satisfies $0<\mhk(r)<1$, and we may conclude that
$\|\mhk\|_q \le \|\mhk\|_2 \le C$ for all $p\in [2,\infty]$.
By the triangle inequality,
\beann
\|  1-\fhk^2\|_2 &\le &
  \| 1-\fhk^2 -\mhk^2 \|_2 + \|\mhk^2\|_2 \\
& \le &  |d| + C,
\eeann
and hence we obtain
$$ \int (1-\fhk)^2 \rdr \le \int (1-\fhk^2)^2\rdr \le C  $$
(since $\fhk \ge 0$.)

Choose a function $\eta\in C^\infty(\RR)$ with
$$  \eta(r) = \cases{  1, & if $r\le 2$, \cr
           0, & if $r\ge 3$, \cr}  $$
and $0\le \eta(r)\le 1$ for all $r$.
Using $\eta^2 \fhk$ as a test function in the weak form
of the rescaled equation for $\fhk$,
\beann
\int \eta^2 \left[  (\fhk')^2 + {d^2\over r^2} \fhk^2\right] \rdr
&=&
\int \left[  (1-\fhk^2 -\mhk^2)\fhk^2\eta^2 - \eta\eta' \fhk\fhk' \right]\rdr \\
&\le &  \int \left[ \frac12 (1-\fhk^2)_2 + \frac12 \eta^4 
          + \frac12 \eta^2(\fhk')^2 + 2 \fhk^2 (\eta')^2 \right]\rdr \\
&\le &  C + \frac12 \int \eta^2 (\fhk')^2\rdr .
\eeann
Absorbing the last term back to the left hand side,
\be\label{L3}
\int \eta^2 \left[  (\fhk')^2 + {d^2\over r^2} \fhk^2\right] \rdr \le C.
\ee

Now choose another smooth function $f_0$ with
$$  f_0(r) = \cases{
   0, & if $r\ge 1$, \cr
   1, & if $r\ge 2$, \cr  }
$$
and $0\le f_0(r) \le 1$.  Note that with this choice
$f_0^2 + \eta^2 \ge 1$.
We use $(\fhk -1)f_0^2$ as a test function in the equation
for $\fhk$ to obtain
\beann
\int \left[  (\fhk')^2 + {d^2\over r^2} (\fhk^2-1)^2 \right]\rdr
   &=& \int \left[ (1-\fhk^2-\mhk^2)\fhk(\fhk-1)f_0^2
       - (\fhk -1)\fhk' f'_0 f_0 \right]\rdr  \\
&& \qquad +\int  \left[- {d^2\over r^2}(\fhk -1)f_0^2
       + {\Shk(2d-\Shk)\over r^2}(\fhk -1)^2 f_0^2 \right]\rdr
\\
&\le &
C\int_1^2 (\fhk-1) |\fhk'| f_0 \rdr  \\
&&\qquad  + 
   \int_1^\infty \left[
       {d^2\over r^2}(1-\fhk)  + {4d\over r^2} (\fhk-1)^2 \right]\rdr \\
&\le & 
\frac12 \int (\fhk')^2f_0^2\rdr 
 + C \|\fk-1\|_2^2 
   + \int \left[ {C\over r^4} + (\fhk^2-1)^4 \right]\rdr \\
&\le & C + \frac12 \int (\fhk')^2f_0^2\rdr .
\eeann
(Note that in the first line, the first integrand is non-positive.)
In conclusion,
\be \label{L4}
\int  \left[  (\fhk')^2 + {d^2\over r^2} (\fhk^2-1)^2 \right]\rdr \le C.
\ee

Now define $\uk = \fhk-f_0\in X$.
Then from (\ref{L3}), (\ref{L4}) we obtain:
\beann
\int \left[  (\uk')^2 + \uk^2 + {d^2\over r^2}\uk^2 \right]\rdr
&\le &
2\int \left[ (\fhk'^2) + (f_0)^2\right]\rdr +
   \int_0^2 {2d^2\over r^2} \fhk^2\rdr \\
&&\qquad 
      + \int_2^\infty (d^2+1)(\fhk-1)^2 \\
&\le &
  2\int \left[ (\eta^2 + f_0^2) (\fhk')^2 
      + {d^2\over r^2}\fhk^2\eta^2 \right] \rdr  + C
\\
&\le & C.
\eeann
In other words, $\uk$ is uniformly bounded in $X$, and
we may extract  weakly convergent subsequences 
$u_n=u_{\kappa_n}\weak u_*$ (in $X$), $m_n=m_{\kappa_n}\to m_*$
(in $H$).  

\noindent
{\bf Step 2:} Strong convergence.
\newline
We next show that the sequences $u_n$, $m_n$ converge in norm.
Let $f_n=f_0+u_n$ and $S_n=\hat S_{\kappa_n}$.  First note that
$$   (1-f_n^2 -m_n^2)m_n  - (1-f_p^2-m_p^2)m_p
   = (1-f_n^2)w - (m_n^2+m_nm_p + m_p^2)w
       + (f_n^2-f_p^2)m_p. $$
Hence, using compact embeddings of $X$, $H$ into $L^q$ for $2<q<\infty$,
\beann
&& \int \left[ ((m_n-m_p)')^2 + g(m_n-m_p)^2\right]\rdr  \\
&&\qquad = \int \left[  
   (1-f_n^2) - (m_n^2+m_nm_p + m_p^2)(m_n-m_p) \right] (m_n-m_p)^2 \rdr \\
&& \qquad\qquad 
 + \int (f_n+f_p)(u_n-u_p)m_p (m_n-m_p) \rdr \\
&&\qquad =  o(1).
\eeann
Therefore, $m_n\to m_*$ in norm.

We proceed in the same way with $u_n$:
\bea \nnn
&&  \int \left\{
       (u'_n-u'_p)^2 
   + \left[ {(d-S_n)^2\over r^2} f_n - {(d-S_p)^2\over r^2} f_p \right](u_n-u_p)
      \right\}\rdr \\
\label{L6}
&& \qquad = 
 \int \left\{ (1-f_n^2-m_n^2)f_n - (1-f_p^2-m_p^2)f_p \right\}(u_n-u_p)\rdr
\eea
Now we expand,
\beann
{(d-S_n)^2\over r^2}f_n - {(d-S_p)^2\over r^2}f_p &=&
         \left[ {d^2\over r^2} - {S_n\over r^2}(2d-S_n)\right](u_n-u_p) \\
&&\qquad + f_p\left[ {S_p\over r^2}(2d-S_p) - {S_n\over r^2} (2d-S_n)\right].
\eeann

Now we take each term separately:
$$
 \int_0^1 {S_n\over r^2}(2d-S_n)(u_n-u_p)^2\rdr
     \le \sup_{r\in [0,1]} |S_n| \int_0^1 2d{(u_n-u_p)^2\over r^2} \rdr \to 0,
$$
since $S_n\to 0$ locally uniformly, and $u_n$ are uniformly bounded.
$$  \int_1^\infty  {S_n\over r^2}(2d-S_n)(u_n-u_p)^2\rdr
  \le \sup \left|{S_n\over r}\right|  \int_1^\infty 2d (u_n-u_p)^2\rdr
   \to 0,
$$
by (\ref{L5}).
Choose $r_0>0$ so that $\int_{r_0}^\infty r^{-2}\, dr < \eps^2/4,$
and $\kappa$ sufficiently large so that 
$$  dr_0^2 \sup \left|{S_p\over r}\right| \|u_n-u_p\|_X < {\eps\over 2}.
$$
Then,
$$
 \int_0^{r_0}
   f_p {S_p\over r^2} (2d-S_p)(u_n-u_p) \rdr
  \le 
dr_0^2 \sup_{0\le r\le r_0} \left|{S_p\over r}\right|
          \sqrt{ \int_0^{r_0} {(u_n-u_p)^2\over r^2}\rdr } \le \eps/2, 
$$
and
$$   \int_{r_0}^\infty  f_p {S_p\over r^2} (2d-S_p)(u_n-u_p)\rdr
    \le 2d^2 \left[ \int_{r_0}^\infty {dr\over r^3} \right]^{1/2}
          \left[ \int_0^\infty (u_n-u_p)^2\rdr \right]^{1/2} <\eps/2.
$$
We return to (\ref{L6}), and substitute the above estimates:
\beann
&& \int \left\{ (u'_n-u'_p)^2 + {d^2\over r^2}(u_n-u_p)^2 \right\}\rdr
   + o(1)
\\
&&\qquad\qquad =
\int \left\{ (1-f_n^2-m_n^2)f_n - (1-f_p^2-m_p^2)f_p \right\}(u_n-u_p)\rdr \\
&& \qquad\qquad =
\int \left\{ (1-3f_0^2) (u_n-u_p)^2 
    -3 f_0 (u_n^2-u_p^2)(u_n-u_p) \right.   \\
&&\qquad\qquad \qquad   \left.
    - (u_n^3-u_p^3)(u_n-u_p)- f_0(m_n^2-m_p^2)(u_n-u_p) \right. \\
&& \qquad\qquad\qquad  \left.
- m_n^2(u_n-u_p)^2 -  u_p(m_n^2-m_p^2)(u_n-u_p) \right\}\rdr \\
&&\qquad\qquad  
= -2\int (u_n-u_p)^2 \rdr + o(1),
\eeann
where we use the facts that $m_n\to m_*$ strongly in $H$,
$u_n$ is bounded in $X$, and $u_n\to u_*$ in $L^2_{loc}$.
In conclusion, the subsequence $u_n\to u_*$ strongly in $X$.

\noindent
{\bf Step 3:} Determining when $m_\infty=0$.
\newline
Since all solutions of \GLhk\ with $g\ge g_\infty^*$
have $m_\infty=0$, we have $\mhk\to 0$ when $g\ge g_\infty^*$.
On the other hand, suppose $\mhk>0$ for all sufficiently large $\kappa$,
but $\mhk\to 0$.  By uniqueness of the normal core solution,
$\fhk\to \finf$, the unique solution of
$$   -\Delta_r \finf + {d^2\over r^2} \finf = (1-\finf^2)\finf. $$
Let
$$    t_\kappa = \int (1-\fhk^2)\mhk^2 \rdr \to 0, $$
and set
$w_\kappa = \mhk/t_\kappa$.  Then
$$   -\Delta_r w_\kappa + g w_\kappa = (1-\fhk^2-\mhk^2)w_\kappa.  $$
Since
$$  \int \left[ (w'_\kappa)^2 + g w_\kappa^2 \right] \rdr
  = \int (1-\fhk^2-\mhk^2)w_\kappa^2 \rdr \le 1  $$
(by the choice of $t_\kappa$,) the bound $\|w_\kappa\|_H\le 1/g$
results.  We extract a subsequence (which we still denote by
$w_\kappa$) with $w_\kappa\weak w_\infty$ weakly in $H$.
Note that $w_\infty\ge 0$.
By the choice of $t_\kappa$, the uniform convergence $\fhk\to\finf$,
and the $L^2_{loc}$ convergence of $w_\kappa\to w_\infty$ we have:
\beann
\int (1-\finf^2)w_\infty^2 \rdr
 &= &  \int \left[
     (1-\finf^2)(w_\infty^2 -w_\kappa^2) + 
           (\fhk^2-\finf^2)w_\kappa^2 
      + (1-\fhk^2) w_\kappa^2  \right]\rdr  \\
&= & 1 + o(1).
\eeann
In particular $w_\infty\not\equiv 0$.
By weak convergence we may pass to the limit in the
equation for $w_\kappa$, and hence $w_\infty$ is a nontrivial
non-negative solution of
$$   -\Delta_r w_\infty + g w_\infty = (1-\finf^2)w_\infty.  $$
This can only occur when $g=g_\infty^*$.

This completes the proof of Theorem~\ref{limit}.
\QED

\section{Estimates and existence}

\setcounter{equation}{0}
\setcounter{th}{0}

In this section we derive the technical estimates which were 
needed in our analysis of the bifurcation problem in 
Section~4.
We also provide the details of the proof of existence of minimizers
of the energies $\EE$ and $\EEhk$.

\subsection{A priori estimates}

We may now prove {\it a priori} estimates for the solutions of our
system \GL, Theorem~\ref{apriori}, as well as the compactness
result for solutions of \GLhk\ (both theorems as stated in the previous section.)
Note that both theorems are stated for all solutions, not only energy minimizers,
and hence we will use our Pohozaev identity (Proposition~\ref{pohozaev})
to obtain energy independent estimates.
As before, we denote by $\fk$, $\Sk$ a normal core solution at $\kappa$,
and $u=f-\fk$, $v=(S-\Sk)/r$.  

By the Pohozaev identity and Lemma 4.2 of \cite{BC} we have
\be\label{firstest}  \kappa^2 \int \left[
           gm^2 + \frac12 (1-f^2-m^2)^2\right]\rdr
        = \int \left( {S'\over r}\right)^2 \rdr \le {d^2\over 2}.  
\ee
In particular, we obtain
$$  \int (1-f)^2 \rdr\le C + C/g, \qquad  \int m^2\rdr \le C/g  $$
with constant $C$ depending on $\kappa, d$.
From the first estimate we  obtain
$$  \|u\|_2 \le \| \fk + u -1\|_2 + \| \fk -1\|_2 \le C+ C/g.  $$
The equation for $m$ together with the second estimate gives:
$$
\int \left[ (m')^2 + \k^2 gm^2 \right]\rdr = \k^2\int (1-f^2-m^2)m^2\rdr 
     \le \kappa^2\int m^2\rdr \le C/g.
$$
In particular, $\|m\|_H \le C$, $\|u\|_2 \le C$, and the constant depending
on $\k$, $d$ may be chosen uniformly
for $g\in J$.

Using the right half of (\ref{firstest}) we have
\be\label{est2}
{d^2\over 2} \ge 
  \int  
    \left({S'\over r}\right)^2 \rdr =
       \int \left[
        \left({\Sk'\over r}\right)^2 + 2 {\Sk'\over r}{(rv)'\over r} +
            \left( {(rv)'\over r}\right)^2\right] \rdr .
\ee
Since
$$
  \left| 2\int {\Sk'\over r}{(rv)'\over r} \rdr \right|
       \le 
          2\int \left({\Sk'\over r}\right)^2 
            + \frac12 \int \left( {(rv)'\over r}\right)^2 \rdr ,
$$
and 
$$  \int \left( {(rv)'\over r}\right)^2 \rdr = 
      \int \left[  (v')^2 + {v^2\over r^2} \right]\rdr,
$$
we may conclude from (\ref{est2}) that
\be\label{est3}
\int \left[  (v')^2 + {v^2\over r^2} \right]\rdr \le C,
\ee
with constant depending only on $d$.  From the embedding
properties of $X$, Lemma \ref{embedding}, we
conclude that $\|v\|_\infty \le C$.

We now use $v$ as a test function in the weak form of the equation for
$S$ to obtain an estimate:
\bea\label{est4}
\left| \int {d-S\over r} f^2 v \rdr\right|
    &=& \left|\int {S'\over r^2} (rv)' \rdr \right| \\
\nnn
&\le &  \frac12 \int \left({S'\over r^2}\right)^2 \rdr
      + \frac12\int \left[  (v')^2 + {v^2\over r^2} \right]\rdr \le C.
\eea
On the other hand, expanding the left-hand side of (\ref{est4}),
\bea\label{est5}
\int \left( {d-S\over r}\right) f^2 v \rdr
    & = &
\int \left[ {d-\Sk\over r} - v\right] (\fk + u)^2 v\rdr \\
\nnn
& = & \int \left( {d-\Sk\over r}\right) [\fk^2 + 2\fk u + u^2] v\rdr
           - \int 2\fk u v^2\rdr\\
\nnn
&&\qquad -\int v^2\fk^2 \rdr -\int v^2 u^2\rdr.
\eea
To bound the term $\int v^2 \fk^2\rdr$, we need to evaluate the
other terms:
\beann
&&  \left| \int {d-\Sk\over r} \fk^2 v \rdr \right|
     \le 2\int \left[ {d-\Sk\over r}\right]^2 \fk^2\rdr
        + \frac18  \int \fk^2 v^2 \rdr
         \le C + \frac18  \int \fk^2 v^2 \rdr,\\
&&   2\left| \int {d-\Sk\over r} \fk u v \rdr \right|
      \le  \int \left[ {d-\Sk\over r}\right]^2 \fk^2\rdr 
            + \|v\|_\infty^2 \|u\|_2^2 \le C, \\
&&  \left| \int {d-\Sk\over r} u^2 v \rdr \right|
     \le \frac12 \|u\|_2^2  \|v\|_\infty^2  +
         \frac12 \int \left[ {d-\Sk\over r}\right]^2 u^2\rdr, \\
&&  2\left| \int \fk u v^2\rdr \right| \le
         8\|u\|_2^2 \|v\|_\infty^2   + \frac18 \int \fk^2 v^2\rdr
        \le C +  \frac18  \int \fk^2 v^2 \rdr,\\
&&\int v^2u^2\rdr\le \|v\|_\infty^2 \|u\|_2^2 \le C.
\eeann
Hence,
\be\label{est6}
\frac34 \int v^2 \fk^2\rdr \le C 
        + \frac12 \int \left[ {d-\Sk\over r}\right]^2 u^2\rdr
\ee

Finally, we use $u$ as a test function in the weak form
of the equation for $f$.  Recalling the definition of $\fk$ as
a normal core solution, we expand and cancel terms to arrive at:
\bea\label{est7}
\int \left[  (u')^2 + \left({d-\Sk\over r}\right)^2\right] \rdr 
& = &\int \left[   2  \left({d-\Sk\over r}\right)(\fk+u) uv
           -  (\fk + u)uv^2 \right] \rdr \\
\nnn
& & +\kappa^2\int \left[
       (1-3\fk^2) u^2 - 3\fk u^3 -u^4 - m^2 f u\right]\rdr.
\eea
Each term on the right hand side may be controlled as
follows:
\beann
&&\left| \int m^2 f u\rdr \right| \le \frac12  \int m^4 + \frac12 \int u^2\rdr \le C,\\
&&\left| \int (1-3\fk^2)u^2 \rdr \right| \le
     3\|u\|_2^2\le C,\\
&&2\left| \int \left({d-\Sk\over r}\right)\fk uv\rdr \right| \le
    \|u\|_2^2 \|v\|_\infty^2 + \int \left({d-\Sk\over r}\right)^2 \fk^2 \rdr \le C,\\
&&\left| \int \fk  u^3 \rdr \right| \le
          \frac32\int \fk^2u^2 \rdr + \frac12 \int u^4\rdr,  \\
&&2\left| \int \left({d-\Sk\over r}\right) u^2 v\rdr \right| \le
    6\|u\|_2^2 \|v\|_\infty^2 + \frac16 \int \left({d-\Sk\over r}\right)^2 u^2 \rdr \\
&&\qquad \le 
        C +  \frac16 \int \left({d-\Sk\over r}\right)^2 u^2 \rdr ,\\
&&\left| \int \fk uv^2\rdr \right| \le \frac12 \|u\|_2^2 \|v\|_\infty^2
           + \frac12 \int \fk^2 v^2\rdr
              \le C +\frac13 \int \left({d-\Sk\over r}\right)^2 u^2 \rdr,
\eeann
where in the last estimate we apply (\ref{est6}).  Using (\ref{est7})
we have
$$  \int \left[    (u')^2 
        + \frac12 \left({d-\Sk\over r}\right)^2 u^2\right] \rdr \le C. $$
Consequently, $\|u\|_X\le C$.
Returning to (\ref{est6}), it follows that
$$ \int \fk^2 v^2\rdr  \le C 
      + \frac23 \int \left({d-\Sk\over r}\right)^2 u^2\rdr  \le C,  $$
and hence (\ref{est3}) yields $\|v\|_X\le C$.
This concludes the proof of Theorem~\ref{apriori}.

\medskip

An analogous result may be proven for solutions of \GLhk:
\begin{th}\label{apriorihk}
Let $d$ be fixed.  For any compact
interval $J\in (0,\infty)$ there exists $C_0=C_0(d,J)>0$ such
that every admissible solution $(f,m)$ of \GLhk\ with $g\in J$
satisfies $\|(f,m)\|_{Z_0}\le C$.
\end{th}
The proof of Theorem~\ref{apriorihk} is similar to (and simpler than) 
the previous one, and
is left to the reader.

\subsection{Compactness}

Here we prove Theorem~\ref{compactness}, which asserts
that the family of solutions to \GLhk\ 
with $g$ bounded away from zero is
a compact set.  The same result holds for \GL, although
the proof is more complicated due to the additional terms
involving $S(r)$.  

\noindent
{\bf Proof of Theorem~\ref{compactness}:}\
Suppose $f_n=\finf + u_n$, $m_n$, $g_n$ are a sequence of
solutions of (GL)${}_{\infty,g_n}$ with $g_n\in [a,b]$.  By the Theorem
we have $\|u_n\|_X, \|m_n\|_H \le C$, and hence we may
extract a subsequence with $u_n\weak \tilde u$, $m_n\weak \tilde m$, and
$g_n\to \tilde g \in [a,b]$.
Then we have
\bea\label{comp1}
&&\int \left[    (u'_n -u'_k)^2 + {d^2\over r^2} (u_n-u_k)^2 \right]\rdr
  =  \\
&&\qquad\int \left[ (1-f_n^2-m_n^2)f_n - (1-f_k^2-m_k^2)f_k \right] (u_n-u_k)\rdr \nonumber \\
\label{comp2}
&&\int \left[    (m'_n -m'_k)^2 + \tilde g (m_n-m_k)^2 \right]\rdr
   =   \\
&&\qquad\int \left[ (1-f_n^2-m_n^2)m_n - (1-f_k^2-m_k^2)m_k \right] (m_n-m_k)\rdr 
            + o(1). \nonumber
\eea
We now expand the two right-hand side terms.
First, we use the embedding properties of $X$, $H$ and the
fact that $0\le f_n <1$ for any solution to show:
\beann
&&\left[ (1-f_n^2-m_n^2)f_n - (1-f_k^2-m_k^2)f_k \right] (u_n-u_k)
= \\
&& \qquad- \int 2 \finf^2(u_n-u_k)^2 \rdr + \int (1-\finf^2)(u_n-u_k)^2 \rdr \\
   &&\qquad          - 2\int [\finf (u_n+u_k)(u_n-u_k)^2
                   -f_n(u_n+u_k)(u_n-u_k)^2]\rdr  \\
&& \qquad    - \int [f_n(m_n+m_k)(m_n-m_k)(u_n-u_k) + (u_k^2+m_k^2)(u_n-u_k)^2]
             \rdr \\
&&= - \int 2 \finf^2(u_n-u_k)^2  \rdr + o(1).
\eeann
Applying the  above estimate to (\ref{comp1}) we have
$$  \int \left[  (u'_n -u'_k)^2 
          + \left({d^2\over r^2}+ 2\finf^2\right)(u_n-u_k)^2 \right]\rdr
     \to 0  $$
as $n,k\to\infty$, so $u_n\to \tilde u$ in  norm on the space $X$.

Similarly, we estimate
\beann
&&\int \left[ (1-f_n^2-m_n^2)m_n - (1-f_k^2-m_k^2)m_k \right] (m_n-m_k)\rdr
          = \\
&& \qquad \int (1-\finf^2)(m_n-m_k)^2 \rdr -\int [ 2\finf u_n (m_n-m_k)^2  \\
&&  \qquad  + 2\finf m_k (m_-m_k)(u_n-u_k) + u_n^2(m_n-m_k)^2]\rdr \\
&& \qquad  -\int [ m_k (u_n^2-u_k^2)(m_n-m_k) 
                   + (m_n^3-m_k^3)(m_n-m_k)]\rdr  \\
&=& o(1).
\eeann
Therefore, (\ref{comp2}) implies that $m_n\to \tilde m$ in $H$.
By passing to the limit in the weak formulation of (GL)${}_{g_n,\infty}$
we easily obtain that $(\tilde f, \tilde m)$ solve (GL)${}_{\tilde g,\infty}$,
and hence the specified solution set is compact.
\QED

\subsection{Existence}

Let $(u_n,v_n,m_n)$ be a minimizing sequence for $\II$, so
$(f_n,S_n,m_n)=(f_0+ u_n,S_0+rv_n,m_n)$ is a minimizing
sequence for $\EE$.  To prove
 Theorem~\ref{existence} we first observe that the energy
$\EE$ is a sum of positive terms, and hence each is individually
bounded.  In particular, $m_n$ is uniformly bounded in $H$.

Now we must estimate $u_n$.  First, note that
$\EE (|f_n|,S_n,m_n) = \EE (f_n,S_nm_n)$, and so we may
assume that our minimizing sequence satisfies $f_n(r)\ge 0$ for all $r$.
Next, we observe 
\be\label{a1}
\|1-f_n^2\|_2 \le \| 1-f_n^2-m_n^2\|_2 + \|m_n^2\|_2 \le
       \| 1-f_n^2-m_n^2\|_2 + C. \ee
Hence we conclude that
\beann C &\ge & \EE(f_n,S_n,m_n) \\
&\ge  &
\int\left\{  (f'_n)^2 + \left[ {S'_n\over r}\right]^2 
     + {(d-S_n)^2 f_n^2\over r^2} + {\kappa^2\over 2}(1-f_n^2)^2 \right\}\rdr.
\eeann
The right-hand side of the above inequality is the free energy of conventional Ginzburg--Landau vortices
studied in \cite{ABG}.  The boundedness of $\|u_n\|_X,\|v_n\|_X$
then follows from the argument of Proposition~4.2 of \cite{ABG}.
We may then pass to the limit in $\EE$ via lower semicontinuity
of the norms and Fatou's Lemma.
\QED

To prove Theorem~\ref{existencehk}
let $(u_n,m_n)$ be a minimizing sequence for $I_\infty$ in $X\times H$,
so $(f_n,m_n)=(\finf+ u_n, m_n)$ is a minimizing sequence
for $\EEhk$.
Choose $r_g\ge 1$ so that ${d^2\over r_g^2}\le {g\over 2}$.  Then
\beann
\EEhk (f_n,m_n) & \ge &
  \int_0^{r_g} \left[  (m_n')^2 + gm_n^2 - {d^2\over r^2}\finf^2 \right]\rdr \\
 &&  \qquad + \int_{r_g}^\infty [
           (m_n')^2 +(g- d^2/r^2) m_n^2 + {d^2\over r^2}(f_n^2+m_n^2-1) \\
 &&  \qquad +\frac12 (f_n^2 +m_n^2 -1)^2 + {d^2\over r^2}(1-\finf^2)] \rdr \\
&\ge & \int_0^\infty  \left[ (m_n')^2 + {g\over 2} m_n^2\right]\rdr
        - \int_0^{r_g} {d^2\over r^2} \finf^2\rdr
           + \int_{r_g}^\infty \left[ {d^2\over r^2}(1-\finf^2) - 
             {d^4\over 2r^4} \right]\rdr,
\eeann
where we have used the elementary bound $ax+x^2/2 \ge -a^2/2$.
In particular, $\EEhk$ is bounded below and the minimizing sequence
has $\|m_n\|_H\le C$ uniformly in $n$.  By the Sobolev embedding,
we also conclude that $\|m_n\|_p\le C_p$ for all $p\in [2,\infty)$.

Now we must estimate $u_n$.  As above we note that
$\EEhk (|f_n|,m_n) = \EEhk (f_n,m_n)$, and so we may
assume that our minimizing sequence satisfies $f_n(r)\ge 0$ for all $r$,
and the bound (\ref{a1}) holds.  Note that we also have:
\be\label{a2}
\|u_n\|_2 \le \|\finf-1\|_2 + \|1-f_n\|_2 \le C + \|1-f_n\|_2.
\ee
By the estimate on $m_n$, (\ref{a1}), and (\ref{a2}) we now have
\beann
C &\ge &  \int \left[  (f'_n)^2 + {d^2\over r^2} (f_n^2- \finf^ 2)
         + \frac12 (1-f_n^2)^2 \right]\rdr \\
 & = & \int \left[ (u'_n)^2 + 2\finf' u'_n + (\finf')^2
       + {d^2\over r^2}(2\finf u_n + u_n^2)
              + \frac12 (1-f_n^2)^2 \right]\rdr \\
& = & \int \left[ (u'_n)^2  + (\finf')^2
       + {d^2\over r^2} u_n^2
              + \frac12 (1-f_n^2)^2 + 2(1-\finf^2)\finf u_n\right]\rdr \\
&\ge & \int \left[  (u'_n)^2 + {d^2\over r^2} u_n^2 + \frac14 u_n^2
     - 8(1-\finf^2)^2 -\frac18 \finf^2 u_n^2 \right]\rdr - C\\
&\ge & \int \left[  (u'_n)^2 + {d^2\over r^2} u_n^2 + \frac18 u_n^2
      \right]\rdr - C.
\eeann
In conclusion $\|u_n\|_X\le C$.  We extract a subsequence for which
both $u_n\weak u_0$ and $m_n\weak m_0$ weakly in $X,H$ respectively,
and pointwise almost everywhere.

By semicontinuity of
the norm, Fatou's Lemma (for the positive terms) and the $L^2_{r,loc}$
convergence of $u_n\to u_0$ we can pass to the limit in (\ref{hkI}):
$$
\IIhk (u_0,m_0) \le \liminf_{n\to\infty} I_\infty (u_n,m_n)
         = \inf_{X\times H} I_\infty.
$$
So the infimum of $I_\infty$ is attained.
\QED

\end{document}